\documentclass[letterpaper, 10 pt, conference]{ieeeconf}  

\IEEEoverridecommandlockouts        
\title{\LARGE \bf Exploring Synchronization in Complex Oscillator Networks}

\author{Florian D\"orfler \and Francesco Bullo%
  \thanks{This material is based in part upon work supported by NSF
    grants IIS-0904501 and CPS-1135819.}%
  \thanks{Florian D\"orfler and Francesco Bullo are with the Center for
    Control, Dynamical Systems and Computation, University of California at
    Santa Barbara. Email: {\tt \{dorfler,bullo\}@engineering.ucsb.edu}} }

\usepackage{graphicx,color}
\usepackage{amsmath}
\usepackage{amssymb}
\usepackage{mathrsfs}
\usepackage[vlined,ruled]{algorithm2e}
\usepackage{subfigure}
\usepackage{cite}

\newtheorem{theorem}{Theorem}[section]
\newtheorem{lemma}[theorem]{Lemma}

\newtheorem{remark}{Remark}

\newcommand*\fvec[1]{\ensuremath{\mathbf{#1}}}                                	
\newcommand*\dt[0]{\frac{d}{d\,t}\,}					                      	
\newcommand{\abs}[1]{\left\lvert{#1}\right\rvert}                     			
\newcommand{\norm}[1]{\left\lVert{#1}\right\rVert}                     		
\newcommand*\mc[0]{\mathcal}                                                  		
\newcommand*\mbb[0]{\mathbb}                                                  		
\DeclareMathOperator*{\sinc}{sinc} 	                                       		
\DeclareMathOperator*{\diag}{\mathrm{diag}}						

\newcommand{\Ker}{\mathrm{Ker\,}}

\DeclareMathOperator*{\sinbf}{\mathrm{\bf sin}} 	                              

\newcommand{\until}[1]{\{1,\dots, #1\}}
\newcommand{\subscr}[2]{#1_{\textup{#2}}}

\newcommand{\setdef}[2]{\{#1 \; | \; #2\}}

\newcommand{\map}[3]{#1: #2 \rightarrow #3}
\newcommand{\union}{\operatorname{\cup}}

\newcommand\oprocendsymbol{\hbox{$\square$}}
\newcommand\oprocend{\relax\ifmmode\else\unskip\hfill\fi\oprocendsymbol}

\newcommand{\mycircle}{\ensuremath{\mathbb S^{1}}}
\newcommand{\sphere}{\ensuremath{\mathbb S}}
\newcommand{\torus}{\ensuremath{\mathbb T}}
\newcommand{\real}{\mathbb{R}}
\newcommand{\complex}{\mathbb{C}}
\newcommand{\rot}{\operatorname{rot}}

\newcommand{\arc}{\mbox{Arc}_{n}}
\newcommand{\bararc}{\overline{\mbox{Arc}}_{n}}
\newcommand{\bal}{\mbox{Bal}_{n}}

\begin{document}
\maketitle
\thispagestyle{empty}
\pagestyle{empty}


\begin{abstract}
The emergence of synchronization in a network of coupled oscillators is a pervasive topic in various scientific disciplines ranging from biology, physics, and chemistry to social networks and engineering applications. A coupled oscillator network is characterized by a population of heterogeneous oscillators and a graph describing the interaction among the oscillators. These two ingredients give rise to a rich dynamic behavior that keeps on fascinating the scientific community.  In this article, we present a tutorial introduction to coupled oscillator networks, we review the vast literature on theory and applications, and we present a collection of different synchronization notions, conditions, and analysis approaches. We focus on the canonical phase oscillator models occurring in countless real-world synchronization phenomena, and present their rich phenomenology. We review a set of applications relevant to control scientists. We explore different approaches to phase and frequency synchronization, and we present a collection of synchronization conditions and performance estimates. For all results we present self-contained proofs that illustrate a sample of different analysis methods in a tutorial style.
\end{abstract}


\section{Introduction}
\label{Section: Introduction}

The scientific interest in synchronization of coupled oscillators can be
traced back to the work by Christiaan Huygens on ``an odd kind sympathy''
between coupled pendulum clocks \cite{CH:1673}, and it still fascinates the
scientific community nowadays \cite{SHS:03,ATW:01}. Within the rich
modeling phenomenology on synchronization among coupled oscillators, we
focus on the canonical model of a continuous-time limit-cycle oscillator
network with continuous and bidirectional coupling.

\begin{figure}[htbp]
	\centering{
	\includegraphics[scale=0.565]{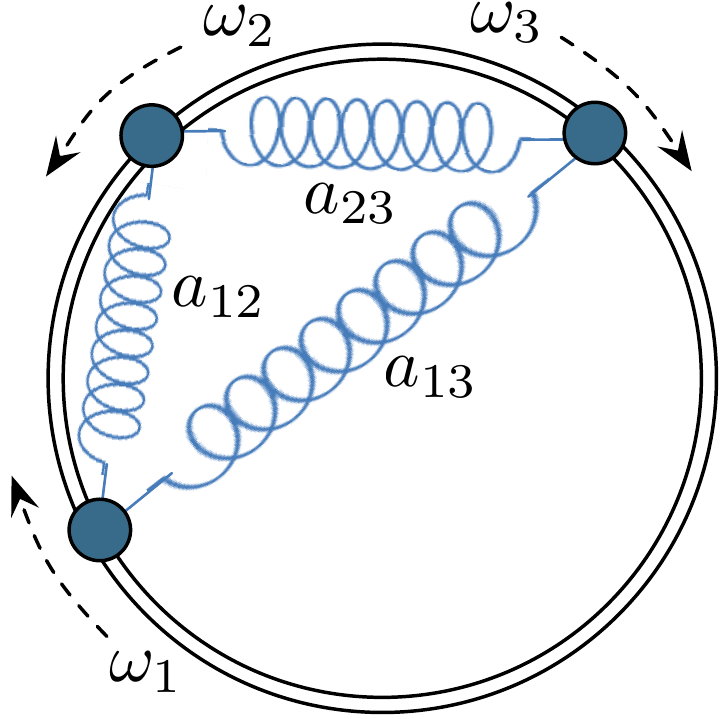}
	\caption{Mechanical analog of a coupled oscillator network}
	\label{Fig: Mechanical analog}
	}
\end{figure}

{\bf A network of coupled phase oscillators:} A mechanical analog of a
coupled oscillator network is the spring network shown in Figure \ref{Fig: Mechanical analog}
and consists of a group of kinematic particles constrained to rotate around
a circle and assumed to move without colliding.
Each particle is characterized by a phase angle $\theta_{i} \in \mathbb
S^{1}$ and has a preferred natural rotation frequency $\omega_{i} \in
\real$. Pairs of interacting particles $i$ and $j$ are coupled through an
elastic spring with stiffness $a_{ij} > 0$. We refer to the Appendix
\ref{Subsection: Modeling of the Mechanical Analog} for a first principle
modeling of the spring-interconnected particles depicted in Figure
\ref{Fig: Mechanical analog}.

Formally, each isolated particle is an oscillator with first-order dynamics
$\dot \theta_{i} = \omega_{i}$.  The interaction among $n$ such oscillators
is modeled by a connected graph $G(\mc V,\mc E,A)$ with nodes $\mc V =
\until n$, edges $\mc E \subset \mc V \times \mc V$, and positive weights
$a_{ij}>0$ for each undirected edge $\{i,j\} \in \mc E$. Under these
assumptions, the overall dynamics of the coupled oscillator network are
\begin{equation}
	\dot \theta_{i}
	=
	\omega_{i} - \sum\nolimits_{j=1}^{n} a_{ij} \sin(\theta_{i}-\theta_{j})
	\,, \qquad i \in \until n
	\label{eq: coupled oscillator model}
	\,.
\end{equation}
The rich dynamic behavior of the coupled oscillator model \eqref{eq:
  coupled oscillator model} arises from a competition between each
oscillator's tendency to align with its natural frequency $\omega_{i}$ and
the synchronization-enforcing coupling $a_{ij} \sin(\theta_{i} -
\theta_{j})$ with its neighbors. Intuitively, a weakly coupled and strongly
heterogeneous network does not display any coherent behavior, whereas a
strongly coupled and sufficiently homogeneous network is amenable to
synchronization, where all frequencies $\dot \theta_{i}(t)$ or even all
phases $\theta_{i}(t)$ become aligned. 
%

{\bf History, applications and related literature:} The coupled oscillator
model \eqref{eq: coupled oscillator model} has first been proposed by
Arthur Winfree \cite{ATW:67}. In the case of a complete interaction graph,
the coupled oscillator dynamics \eqref{eq: coupled oscillator model} are
nowadays known as the {\em Kuramoto model} of coupled oscillators due to Yoshiki
Kuramoto \cite{YK:75,YK:84}. Stephen Strogatz provides an excellent
historical account in \cite{SHS:00}. We also recommend the survey~\cite{JAA-LLB-CJPV-FR-RS:05}.

Despite its apparent simplicity, the coupled oscillator model \eqref{eq:
  coupled oscillator model} gives rise to rich dynamic behavior. This model
is encountered in various scientific disciplines ranging from natural
sciences over engineering applications to social networks.
%
The model and its variations appear in the study if biological
synchronization phenomena such as pacemaker cells in the heart
\cite{DCM-EPM-JJ:87}, circadian rhythms \cite{CL-DRW-SHS-RSM:97},
neuroscience \cite{FV-JPL-ER-JM:01,EB-PH-JM:03,SMC-GBE-MCV-JMB:97}, metabolic synchrony in
yeast cell populations \cite{AKG-BC-EKP:71}, flashing fireflies
\cite{JB:88}, chirping crickets \cite{TJW:69}, biological locomotion \cite{NK-GBE:88}, animal flocking behavior \cite{NEL-TS-NB-LS-IDC-SAL:12}, fish schools \cite{DAP-NEL-RS-DG-JKP:07}, and rhythmic applause
\cite{ZN-ER-TV-YB-AIB:00}, among others.
The coupled oscillator model \eqref{eq: coupled oscillator model} also
appears in physics and chemistry in modeling and analysis of spin glass
models \cite{HD:92,GJ-JA-DB-ACCC-CPV:01}, flavor evolutions of neutrinos
\cite{JP:98}, coupled
Josephson junctions \cite{KW-PC-SHS:98}, and in the analysis of chemical oscillations
\cite{IZK-YZ-JLH:02}.

Some technological applications of the coupled oscillator model \eqref{eq:
  coupled oscillator model} include deep brain stimulation \cite{PAT:03,AN-JM:11},
vehicle coordination
\cite{DAP-NEL-RS-DG-JKP:07,RS-DP-NEL:07,RS-DP-NEL:08,DJK:08,DJK-PL-KAM-TJ:08}, 
carrier synchronization without phase-locked loops \cite{RMMU-RM-SD:11}, semiconductor
lasers \cite{GK-AGV-PM:00,FCH-EMI:00}, microwave oscillators \cite{RAY-RCC:02}, clock
synchronization in decentralized computing networks \cite{WCL-FG-WCH-KD:85,OS-US-YBN-SS:08,YWH-AS:05,RB-AC-LQ-SS-STP:10,YW-NF-JFD:12,EM-AT:11b}, decentralized maximum likelihood estimation \cite{SB-SG:07}, and
droop-controlled inverters in microgrids \cite{JWSP-FD-FB:12j}.
Finally, the coupled oscillator model \eqref{eq: coupled oscillator model}
also serves as the prototypical example for synchronization in complex
networks \cite{AA-ADG-JK-YM-CZ:08,SB-VL-YM-MC-DUH:06,SHS:01,JAKS-GVO:08} and its
linearization is the well-known consensus protocol studied in networked
control, see the surveys and monographs
\cite{ROS-JAF-RMM:07,WR-RWB-EMA:07,FB-JC-SM:09}.
Various control scientists explored the coupled oscillator model \eqref{eq: coupled oscillator model} as a nonlinear generalization of the consensus protocol \cite{LM:05,AJ-NM-MB:04,NC-MWS:08,ZL-BF-MM:07,AS-RS:07,LS-AS-RS:06,ROS:06b}.

Second-order variations of the coupled oscillator model \eqref{eq: coupled
  oscillator model} appear in synchronization phenomena, in population of flashing
fireflies \cite{GBE:91}, in particle models mimicking animal flocking
behavior \cite{SYH-EJ-MJK:10,SYH-CL-BR-MS:11}, in structure-preserving power system models,
\cite{ARB-DJH:81,PWS-MAP:98} in network-reduced power system models
\cite{HDC-CCC-GC:95,FD-FB:09z}, in coupled metronomes \cite{JP:02}, 
in pedestrian crowd synchrony on London's Millennium bridge \cite{SHS-DMA-AMR-BE-EO:05}, and in
Huygen's pendulum coupled clocks \cite{MB-MFS-HR-KW:02}. Coupled oscillator networks with second-order dynamics have been theoretically analyzed in
\cite{YPC-SYH-SBH:11,HAT-AJL-SO:97,HAT-AJL-SO:97b,HH-GSJ-MYC:02,HH-MYC-JY-KSS:99,JAA-LLB-RS:00,JAA-LLB-CJPV-FR-RS:05,FD-FB:10w},
among others.
 
Coupled oscillator models of the form \eqref{eq: coupled oscillator model}
are also studied from a purely theoretic perspective in the physics,
dynamical systems, and control communities. At the heart of the coupled
oscillator dynamics is the {transition from incoherence to
  synchrony}. Here, different notions and degrees of synchronization can be
distinguished \cite{FD-FB:10w,FDS-DA:07,RM-SHS:07}, and the (apparently)
incoherent state features rich and largely unexplored dynamics as well
\cite{YLM-OVP-PAT:05,RT:07,OVP-YLM-PAT:05,JAKS-GVO:08}. In this article we will be
particularly interested in phase and frequency synchronization when all
phases $\theta_{i}(t)$ become aligned, respectively all frequencies $\dot
\theta_{i}(t)$ become aligned. We refer to
\cite{RM-SHS:07,FD-FB:09z,DAP-NEL-RS-DG-JKP:07,JL:11,DJK-PL-KAM-TJ:08,EC-PM:08,MV-OM:09,NC-MWS:08,SYH-TH-KJH:10,AJ-NM-MB:04,REM-SHS:05,DA-JAR:04,MV-OM:08,GBE:85,VM-OM:11,JCB-LDV-MJP:11,JAR-DA:04,FDS-DA:07,SHS:00,JAA-LLB-CJPV-FR-RS:05,FD-FB:10w,LS-AS-RS:06,LDV:11,UM-AP-FA:09,RS-DP-NEL:07,EM-AT:10,LS:10,AF-AC-WPL:11,SYH-MS:11,LB-LS-ADG:09,YM-AFP:04,ACK:10,AS:09,EAC-PAM-RF:10-a,EAC-PAM-RF:10-b,PM:06,JLvH-WFW:93,SJC-JJS:10b,GSS-UM-FA:09,EC-PM:09,SHS-REM:88,AF-AC-WPL:11,CH:10,EAM-EB-SHS-PS-TMA:09,YH-PGM-SPM-UVS:12,DAW-SHS-MG:06,YW-JFD:11b,FD-MC-FB:11v-arxiv}
for an incomplete overview concerning numerous recent research activities.
We will review some of literature throughout the paper and refer to the
surveys
\cite{JAA-LLB-CJPV-FR-RS:05,FD-FB:10w,SHS:00,AA-ADG-JK-YM-CZ:08,SB-VL-YM-MC-DUH:06,SHS:01}
for further applications and numerous additional theoretic results
concerning the coupled oscillator model \eqref{eq: coupled oscillator
  model}.

{\bf Contributions and contents:} In this paper, we introduce the reader to
synchronization in networks of coupled oscillators. We present a sample of
important analysis concepts in a tutorial style and from a
control-theoretic perspective. 

In Section \ref{Section: Technological Applications of Kuramoto Oscillators}, we will review a set of selected technological applications which are directly tied to the coupled oscillator model \eqref{eq: coupled oscillator model} and also relevant to control systems. We will cover vehicle coordination, electric power networks, and clock synchronization in depth, and also justify the importance of the coupled oscillator model \eqref{eq: coupled oscillator model} as a canonical model.
Prompted by these applications, we review the existing results concerning phase
synchronization, phase balancing, and frequency synchronization, and we also present some novel results on
synchronization in sparsely-coupled networks.

In particular, Section \ref{Section: Models and Synchronization Notions}
introduces the reader to different synchronization notions, performance
metrics, and synchronization conditions. We illustrate these results with a
simple yet rich example that nicely explains the basic phenomenology in
coupled oscillator networks.

Section \ref{Section: Analysis of Synchronization} presents a collection of
important results regarding phase synchronization, phase balancing, and frequency synchronization. By now the
analysis methods for synchronization have reached a mature level, and we
present simple and self-contained proofs using a sample of different
analysis methods.
In particular, we present one result on phase synchronization and one result on phase balancing including
estimates on the exponential synchronization rate and the region of
attraction (see Theorem \ref{Theorem: Phase sync} and Theorem \ref{Theorem: Phase balancing}). We also present some
implicit and explicit, and necessary and sufficient conditions for
frequency synchronization in the classic homogeneous case of a complete and
uniformly-weighted coupling graphs (see Theorem \ref{Theorem: Kuramoto
  frequency sync}). Concerning frequency synchronization in sparse graphs, we present two partially new synchronization conditions
depending on the algebraic connectivity (see Theorem \ref{Theorem:
  Synchronization in coupled oscillators I} and Theorem \ref{Theorem:
  Synchronization in coupled oscillators II}).

In our technical presentation, we try to strike a balance between mathematical precision and removing unnecessary technicalities. For this reason some proofs are reported in the appendix and others are only sketched here with references to the detailed proofs elsewhere. Hence, the main technical ideas are conveyed while the tutorial value is maintained.

Finally, Section \ref{Section: Conclusions} concludes the paper. We summarize the limitations of existing analysis methods and
suggest some important directions for future research.

{\bf Preliminaries and notation:} The remainder of this section introduces
some notation and recalls some preliminaries.

{\it Vectors and functions:} Let $\fvec 1_{n}$ and $\fvec 0_{n}$ be the
$n$-dimensional vector of unit and zero entries, and let $\fvec 1_n^\perp$
be the orthogonal complement of $\fvec 1_{n}$ in $\real^{n}$, that is,
$\fvec 1_n^\perp \triangleq \{ x \in \real^{n} :\, x \perp \fvec 1_{n}\}$.
Given an $n$-tuple $(x_{1},\dots,x_{n})$, let $x \in \real^{n}$ be the
associated vector with maximum and minimum elements $\subscr{x}{max}$ and
$\subscr{x}{min}$. For an ordered index set $\mc I$ of cardinality $|\mc
I|$ and an one-dimensional array $\{x_{i}\}_{i \in \mc I}$, let
$\diag(\{c_{i}\}_{i \in \mc I}) \in \real^{|\mc I| \times |\mc I|}$ be the
associated diagonal matrix.
Finally, define the continuous function  $\map{\sinc}{\real}{\real}$ by
$\sinc(x) = \sin(x)/x$ for $x\neq0$.

{\it Geometry on the $n$-torus:} The set $\sphere^{1}$ denotes the {\em unit
  circle}, an {\it angle} is a point $\theta \in \sphere^{1}$, and an {\it
  arc} is a connected subset of $\sphere^{1}$.
The {\em geodesic distance} between two angles $\theta_{1}$, $\theta_{2}
\in \sphere^{1}$ is the minimum of the counter-clockwise and the clockwise
arc lengths connecting $\theta_{1}$ and $\theta_{2}$.  With slight abuse of
notation, let $|\theta_{1}-\theta_{2}|$ denote the {\it geodesic distance}
between two angles $\theta_{1},\theta_{2} \in \sphere^{1}$.
The {\em $n$-torus} is the product set $\torus^{n} = \sphere^{1} \times
\dots \times \sphere^{1}$ is the direct sum of $n$ unit circles.
%
For $\gamma\in{[0,2\pi[}$, let $\bararc(\gamma) \subset \torus^{n}$
be the closed set of angle arrays $\theta = (\theta_1,\dots,\theta_n)$ with
the property that there exists an arc of length $\gamma$ containing all
$\theta_1,\dots,\theta_n$. Thus, an angle array $\theta \in
\bararc(\gamma)$ satisfies $\max\nolimits_{i,j \in \until n} |\theta_{i}
- \theta_{j}| \leq \gamma$. Finally, let $\arc(\gamma)$ be the interior
of the set $\bararc(\gamma)$.

{\it Algebraic graph theory:} Let $G(\mc V,\mc E,A)$ be an undirected,
connected, and weighted graph without self-loops. Let $A \in \mathbb R^{n
  \times n}$ be its symmetric nonnegative {\em adjacency matrix} with zero
diagonal, $a_{ii}=0$.  For each node $i \in \until n$, define the nodal
degree by $\textup{deg}_{i} = \sum\nolimits_{j=1}^{n} a_{ij}$. Let $L \in
\real^{n \times n}$ be the {\em Laplacian matrix} defined by $L =
\diag(\{\textup{deg}_{i} \}_{i =1}^{n}) - A$.
If a number $\ell \in \until {|\mc E|}$ and an arbitrary direction is
assigned to each edge $\{i,j\} \in \mc E$, the (oriented) {\em incidence
  matrix} $B\in \real^{n \times |\mc E|}$ is defined component-wise by
$B_{k\ell} = 1$ if node $k$ is the sink node of edge ${\ell}$ and by
$B_{k\ell} = -1$ if node $k$ is the source node of edge ${\ell}$; all other
elements are zero.  For $x \in \real^{n}$, the vector $B^{T}x$ has
components $x_{i} - x_{j}$ corresponding to the oriented edge from $j$ to
$i$, that is, $B^{T}$ maps node variables $x_{i}$, $x_{j}$ to incremental
edge variables $x_{i}-x_{j}$.  If $\diag(\{a_{ij}\}_{\{i,j\} \in \mc E})$
is the diagonal matrix of edge weights, then $L = B \diag(\{a_{ij}
\}_{\{i,j\} \in \mc E}) B^{T}$.
%
%
If the graph is connected, then $\Ker(B^{T}) = \Ker(L) =
\mathrm{span}(\fvec 1_{n})$, all $n-1$ non-zero eigenvalues of $L$ are
strictly positive, and the second-smallest eigenvalue $\lambda_{2}(L)$ is
called the {\it algebraic connectivity} and is a spectral connectivity
measure.


\section{Applications of Kuramoto Oscillators Relevant to Control Systems}
\label{Section: Technological Applications of Kuramoto Oscillators}

The mechanical analog in Figure~\ref{Fig: Mechanical analog} 
provides an intuitive illustration of the coupled oscillator dynamics \eqref{eq: coupled oscillator model}, and we reviewed a wide range of examples from physics, life sciences, and technology in Section \ref{Section: Introduction}. Here, we detail a set of selected technological applications which are relevant to control systems scientists.

\subsection{Flocking, Schooling, and Planar Vehicle Coordination}

An emerging research field in control is the coordination of autonomous vehicles based on locally available information and inspired by biological flocking phenomena. Consider a set of $n$ particles in the plane $\real^{2}$, which we identify with the complex plane $\complex$. Each particle $i \in \mc V = \until n$ is characterized by its position $r_{i} \in \complex$, its heading angle $\theta_{i} \in \mbb S^{1}$, and a steering control law $u_{i}(r,\theta)$ depending on the position and heading of itself and other vehicles. For simplicity, we assume that all particles have constant and unit speed. The particle kinematics are then given by \cite{EWJ-PSK:04}
\begin{align}
	\begin{split}
	\dot r_{i} &= e^{\textup{i} \theta_{i}} \,,\\
	\dot \theta_{i} &= u_{i}(r,\theta) \,,
	\end{split}
	 \qquad \Biggr\} \quad i \in \until n \,,
	 \label{eq: vehicle coordination dynamics}
\end{align}
where $\textup{i} = \sqrt{-1}$ is the imaginary unit. If the control $u_{i}$ is identically zero, then particle $i$ travels in a straight line with orientation $\theta_{i}(0)$, and if $u_{i} = \omega_{i} \in \real$ is a nonzero constant, then the particle traverses a circle with radius $1/|\omega_{i}|$. 

\begin{figure}[htbp]
	\centering{
	\includegraphics[width=0.98\columnwidth]{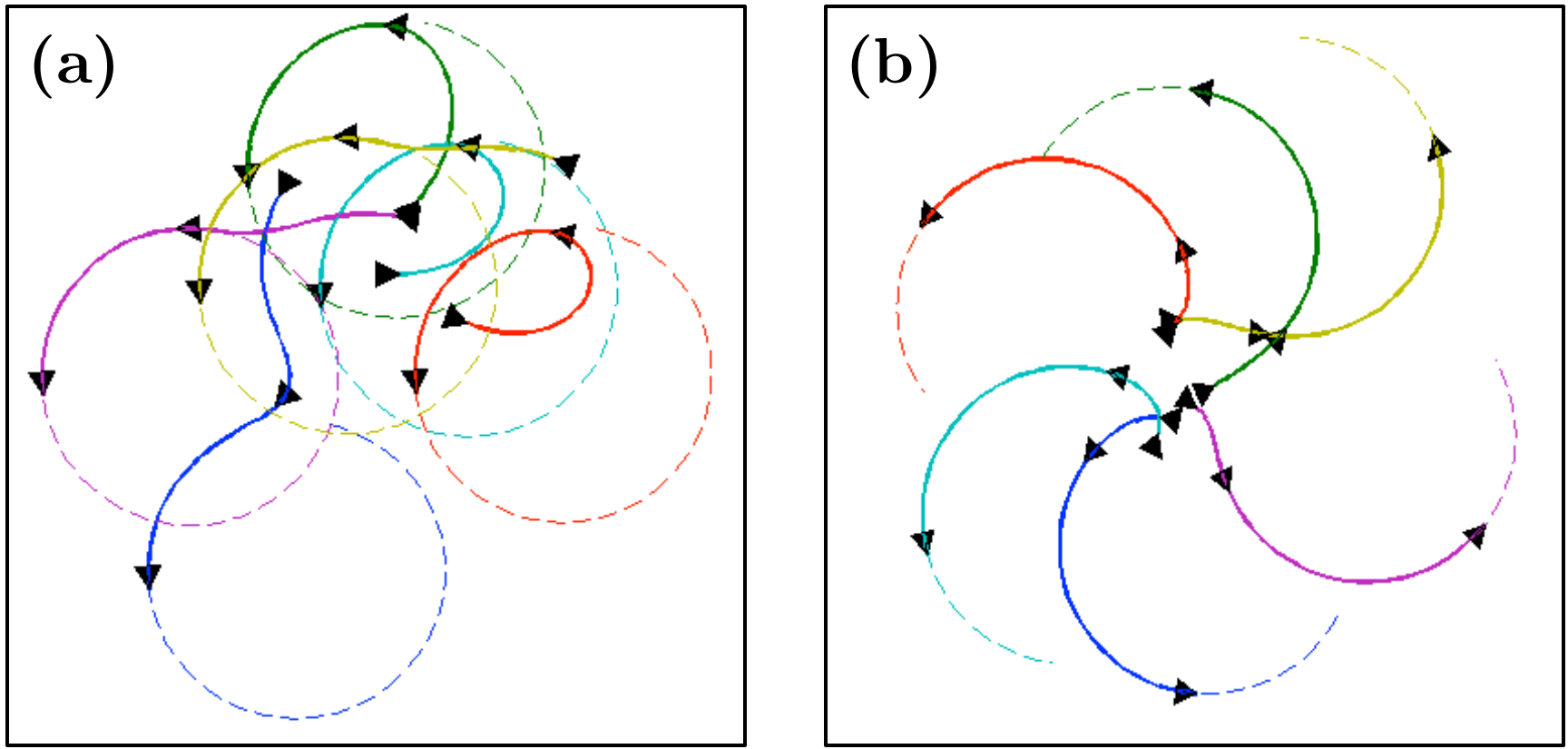}
	\caption{Illustration of the controlled dynamics \eqref{eq: vehicle coordination dynamics}-\eqref{eq: coordination law} with $n\!=\!6$~particles, a complete interaction graph, and identical and constant natural frequencies $\omega_{0}(t) = 1$, where $K\!=\!1$ in panel (a) and  $K\!=\!-1$ in panel~(b). The arrows depict the orientation, the dashed curves show the long-term position dynamics, and the solid curves show the initial transient position dynamics. It can be seen that even for this simple choice of controller, the resulting motion results in ``synchronized'' or ``balanced'' heading angles for $K = \pm 1$.}
	\label{Fig: Vehicle coordination}
	}
\end{figure}

The interaction among the particles is modeled by a possibly time-varying interaction graph $G(\mc V,\mc E(t),A(t))$ determined by communication and sensing patterns. Some interesting motion patterns emerge if the controllers use only relative phase information between neighboring particles, that is, $u_{i} = \omega_{0}(t) + f_{i}(\theta_{i}-\theta_{j})$ for $\{i,j\} \in \mc E(t)$ and $\omega_{0}:\, \real_{\geq 0} \to \real$. For example, the control $u_{i} = \omega_{0}(t) - K \cdot \sum_{j=1}^{n} a_{ij}(t) \sin(\theta_{i} - \theta_{j})$ with gain $K \in \real$ results in 
\begin{equation}
	\dot \theta_{i}
	=
	\omega_{0}(t) - K \cdot \sum\nolimits_{j=1}^{n} a_{ij}(t) \sin(\theta_{i} - \theta_{j})
	\,,\quad i \in \mc V \,.
	\label{eq: coordination law}
\end{equation}
The controlled phase dynamics \eqref{eq: coordination law} correspond to the coupled oscillator model \eqref{eq: coupled oscillator model} with a time-varying interaction graph with weights $K \cdot a_{ij}(t)$ and identically time-varying natural frequencies $\omega_{i} = \omega_{0}(t)$ for all $i \in \until n$. 
The controlled phase dynamics \eqref{eq: coordination law} give rise to very interesting coordination patterns that mimic animal flocking behavior \cite{NEL-TS-NB-LS-IDC-SAL:12} and fish schools \cite{DAP-NEL-RS-DG-JKP:07}. Inspired by these biological phenomena, the controlled phase dynamics \eqref{eq: coordination law} and its variations have also been studied in the context of tracking and formation controllers in swarms of autonomous vehicles \cite{DAP-NEL-RS-DG-JKP:07,RS-DP-NEL:07,RS-DP-NEL:08,DJK:08,DJK-PL-KAM-TJ:08}.
A few trajectories are illustrated in Figure~\ref{Fig: Vehicle coordination}, and we refer to \cite{DAP-NEL-RS-DG-JKP:07,RS-DP-NEL:07,RS-DP-NEL:08,DJK:08,DJK-PL-KAM-TJ:08} for other control laws and motion patterns.

In the following sections, we will present various tools to analyze the motion patterns in Figure~\ref{Fig: Vehicle coordination}, which we will refer to as {\em phase synchronization} and {\em phase balancing}.

\subsection{Power Grids with Synchronous Generators and Inverters}

Here, we present the {\em structure-preserving power network model} introduced in \cite{ARB-DJH:81} and refer to \cite[Chapter 7]{PWS-MAP:98} for detailed derivation from a higher order first principle model. Additionally, we equip the power network model with a set of inverters and refer to \cite{JWSP-FD-FB:12j} for a detailed~modeling.

Consider an alternating current (AC) power network modeled as an undirected, connected, and weighted graph with node set $\mc V = \until n$, transmission lines $\mc E \subset \mc V \times \mc V$, and admittance matrix $Y\!=\!Y^{T} \in \complex^{n \times n}$. For each node, consider the voltage phasor $V_{i} = |V_{i}| e^{\textup{i} \theta_{i}}$ corresponding to the phase $\theta_{i} \in \mbb S^{1}$ and magnitude $|V_{i}| \geq 0$ of the sinusoidal solution to the circuit equations. If the network~is~lossless, then the active power flow from node $i$ to $j$ is $a_{ij} \sin(\theta_{i} - \theta_{j})$, where we used the shorthand $a_{ij} = |V_{i}| \cdot |V_{j}| \cdot \Im(Y_{ij})$.

In the following, we assume that the node set is partitioned as $\mc V = \mc V_{1} \union \mc V_{2} \union \mc V_{3}$, where $\mc V_{1}$ are load buses, $\mc V_{2}$ are conventional synchronous generators, and $\mc V_{3}$ are grid-connected direct current (DC) power sources, such as solar farms. 
The active power drawn by a load $i \in \mc V_{1}$ consists of a constant term $P_{\textup{l},i} > 0$
and a frequency-dependent term $D_{i} \dot \theta_{i}$ with $D_{i}>0$. The
resulting power balance equation is
\begin{equation}
	D_{i} \dot \theta_{i} + P_{\textup{l},i}
	=
	- \sum\nolimits_{j=1}^{n} a_{ij} \sin(\theta_{i} - \theta_{j})
	\,,\;\;\; i \in \mc V_{1} \,.
	\label{eq: load dynamics}
\end{equation}
If the generator reactances are absorbed into the admittance matrix, then the swing dynamics of generator $i \in \mc V_{2}$ are
\begin{equation}
	M_{i} \ddot \theta_{i} + D_{i} \dot \theta_{i}
	= 
	P_{\textup{m},i} - \sum\nolimits_{j=1}^{n} a_{ij} \sin(\theta_{i} - \theta_{j})
	\,,\;\;\;  i \in \mc V_{2} ,
	\label{eq: generator dynamics}
\end{equation}
where $\theta_{i} \in \mbb S^{1}$ and $\dot \theta_{i} \in \real^{1}$ are the generator rotor angle and frequency, $P_{\textup{m},i} > 0$ is the mechanical power input, and $M_{i} > 0$, and $D_{i} > 0$ are the inertia and damping coefficients. 

\begin{figure}[htbp]
	\centering{
	\includegraphics[width=0.98\columnwidth]{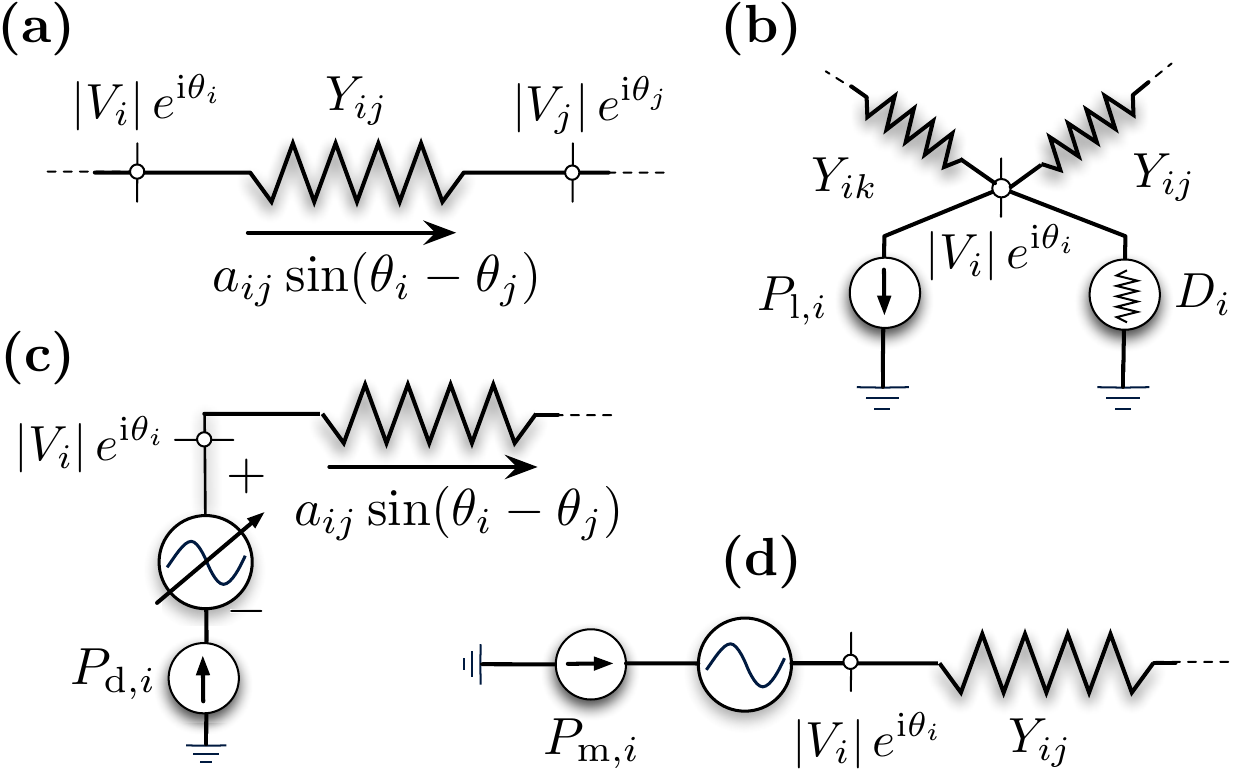}
	\caption{Illustration of the power network devices as circuit elements.~Subfigure (a) shows a transmission element connecting nodes $i$ and $j$, Subfigure (b) shows a frequency-dependent load, Subfigure (c) shows an inverter controlled according to \eqref{eq: inverter dynamics}, and Subfigure (d) shows a synchronous generator.}
	\label{Fig: circuit elements}
	}
\end{figure}

We assume that each DC source is connected to the AC grid via an DC/AC inverter, the inverter output impendances are absorbed into the admittance matrix, and each inverter is equipped with a conventional droop-controller. For a droop-controlled inverter $i \in \mc V_{3}$ with droop-slope $1/D_{i} > 0$, the deviation of the power output $\sum\nolimits_{j=1}^{n} a_{ij} \sin(\theta_{i} - \theta_{j})$ from its nominal value $P_{\textup{d},i} > 0$ is proportional to the frequency deviation $D_{i} \dot \theta_{i}$. This gives rise to the inverter dynamics
\begin{equation}
	D_{i} \dot \theta_{i}
	=
	P_{\textup{d},i}	- \sum\nolimits_{j=1}^{n} a_{ij} \sin(\theta_{i} - \theta_{j})
	\,,\;\;\; i \in \mc V_{3} \,.
	\label{eq: inverter dynamics}
\end{equation}
These power network devices are illustrated in Figure~\ref{Fig: circuit elements}.
Finally, we remark that different load models such as~constant power/current/susceptance loads and synchronous motor loads can be modeled and analyzed by the same set of equations \eqref{eq: load dynamics}-\eqref{eq: inverter dynamics}, see \cite{PWS-MAP:98,SS-PV:80,FD-FB:11d,HDC-CCC-GC:95,FD-FB:09z}. 

Synchronization is pervasive in the operation of power networks. All generating units of an interconnected grid must remain in strict {frequency synchronism} while continuously following demand and rejecting disturbances. 
Notice that, with exception of the inertial terms $M_{i} \ddot \theta_{i}$ and the possibly non-unit coefficients $D_{i}$, the power network dynamics \eqref{eq: load dynamics}-\eqref{eq: inverter dynamics} are a perfect electrical analog of the coupled oscillator model \eqref{eq: coupled oscillator model} with $\omega = (-P_{\textup{l},i},P_{\textup{m},i},P_{\textup{d},i})$. 
Thus, it is not surprising that scientists from different disciplines recently advocated coupled oscillator approaches to analyze synchronization in power networks \cite{VF-SR-EC-EM-VR:09,GF-AHN-NFP:08,JWSP-FD-FB:12j,MR-AS-MT-DW:12,FD-MC-FB:11v-arxiv,FD-FB:09z,DS-UR-BS-MP:01,DJH-GC:06,LB-LS-ADG:09,HAT-AJL-SO:97}. 

The theoretic tools presented in the following sections establish how {\em frequency synchronization} in power networks depend on the nodal parameters $(P_{\textup{l},i},P_{\textup{m},i},P_{\textup{d},i})$ as well as the interconnecting electrical network with weights~$a_{ij}$. Ultimately, this deep understanding of synchrony gives us the correct intuition to design controllers and remedial action schemes preventing the loss of synchrony.

\subsection{Clock Synchronization in Decentralized Networks}

Another emerging technological application of the coupled oscillator model \eqref{eq: coupled oscillator model} is clock synchronization in decentralized computing networks, such as wireless and distributed software networks. A natural approach to clock synchronization is to treat each clock as a coupled oscillator and follow a diffusion-based protocol to synchronize them, see the historic and recent surveys \cite{WCL-FG-WCH-KD:85,OS-US-YBN-SS:08}, the landmark paper \cite{YWH-AS:05}, and the interesting recent results \cite{RB-AC-LQ-SS-STP:10,YW-NF-JFD:12,EM-AT:11b}.
  
Consider a set of distributed processors $\mc V = \until n$ interconnected in a (possibly directed) communication network. Each processor is equipped with an internal software clock, and these clocks need to be synchronized for distributed computing and network routing tasks.
For simplicity, we consider only analog clocks with continuous coupling since digital clocks are essentially discretized analog clocks and pulse-coupled clocks can be modeled continuously after a phase reduction and averaging analysis. 

For our purposes, the clock of processor $i$ is a voltage-controlled oscillator which outputs a harmonic waveform $s_{i}(t) = \sin(\theta_{i}(t))$, where $\theta_{i}(t)$ is the accumulated instantaneous phase. For uncoupled nodes, the phase $\theta_{i}(t)$ evolves~as
\begin{equation*}
	\theta_{i}(t) = \left( \theta_{i}(0) + \frac{2 \pi}{\subscr{T}{nom} + T_{i}} t  \right) \textup{mod}(2\pi)
	\,,\quad i \in \until n \,.
\end{equation*}
where $\subscr{T}{nom} > 0$ is the nominal period, $T_{i} \in \real$ is an offset (frequency offset or skew), and $\theta_{i}(0) \in \mathbb S^{1}$ is the initial phase. 
To synchronize their internal clocks, the processors follow a diffusion-based protocol. In a first step, neighboring oscillators continuously communicate their respective waveforms $s_{i}(t)$ to another. Second, through a phase detector each node measures a convex combination of phase differences as
\begin{equation*}
	\textup{cvx}_{i}(\theta(t)) = \sum\nolimits_{j=1}^{n} a_{ij} f(\theta_{i}(t) - \theta_{j}(t))
	\,,\quad i \in \until n \,,
\end{equation*}
where $a_{ij} \!\geq\! 0$ are convex ($\sum_{j=1}^{n} a_{ij} \!=\! 1$) and detector-specific weights, and $f: \mathbb S^{1} \to \real$ is an odd $2\pi$-periodic~function. Finally, $\textup{cvx}_{i}(\theta(t))$ is fed to a (first-order and constant) phase-locked loop filter $K$ whose output drives the local phase according to
 \begin{equation}
 	\dot \theta_{i}(t) = \frac{2 \pi}{T_{i}} + K \cdot \textup{cvx}_{i}(\theta(t))
	\,, \quad i \in \until n \,.
	\label{eq: sync protocol}
 \end{equation}
 The goal of the synchronization protocol \eqref{eq: sync protocol} is to synchronize the frequencies $\dot \theta_{i}(t)$ or even the phases $\theta_{i}(t)$ in the processor network. 
 For an undirected communication protocol, symmetric weights $a_{ij} = a_{ji}$, and a sinusoidal coupling function $f(\cdot) = \sin(\cdot)$, the synchronization protocol \eqref{eq: sync protocol} equals again the coupled oscillator model \eqref{eq: coupled oscillator model}. 
 
 The tools developed in the next section will enable us to state conditions when the protocol \eqref{eq: sync protocol} successfully achieves phase or frequency synchronization. Of course, the protocol \eqref{eq: sync protocol} is merely a starting point, more sophisticated phase-locked loop filters can be constructed to enhance steady-state deviations from synchrony, and communication and phase noise as well as time-delays can be considered in the~design.

\subsection{Canonical Coupled Oscillator Model}

The importance of the coupled oscillator model \eqref{eq: coupled oscillator model} does not stem only from the various examples listed in Sections \ref{Section: Introduction} and \ref{Section: Technological Applications of Kuramoto Oscillators}. Even though  model \eqref{eq: coupled oscillator model} appears to be quite specific (a phase oscillator with constant driving term and continuous, diffusive, and sinusoidal coupling), it is the {\em canonical model} of coupled limit-cycle oscillators \cite{FCH-EMI:97}.
In the following, we briefly sketch how such general models can be reduced to model \eqref{eq: coupled oscillator model}. We schematically follow the approaches \cite[Chapter 10]{EMI:07},\cite{EMI-YK:06} developed in the computational neuroscience community without aiming at mathematical precision, and we refer to \cite{FCH-EMI:97,GBE-NK:84} for further details.

Consider an oscillator modeled as a dynamical system with state $x \in \real^{m}$ and nonlinear dynamics $\dot x = f(x)$, which admit a locally exponentially stable periodic orbit $\gamma \subset \real^{m}$ with period $T>0$. By a change of variables, any trajectory in a local neighborhood of $\gamma$ can be characterized by a phase variable $\varphi \in \mathbb S^{1}$ with dynamics $\dot \varphi = \Omega$, where $\Omega = 2\pi/T$. 

Now consider a weakly forced oscillator of the form
\begin{equation}
	\dot x = f(x) + \epsilon \cdot \delta(t)
	\label{eq: weakly forced oscillator}
	\,,
\end{equation}
where $\epsilon > 0$ is sufficiently small and $\delta(t)$ is a time-dependent forcing term. For small forcing $\epsilon \delta(t)$, the attractive limit cycle $\gamma$ persists, and the phase dynamics are obtained~as
\begin{equation*}
	\dot \varphi = \Omega + \epsilon Q(\varphi) \delta(t) + \mc O(\epsilon^{2})
	\,,
\end{equation*}
where $Q(\varphi)$ is the infinitesimal phase response curve (or linear response function), and we dropped higher order terms. 

Now consider $n$ such limit cycle oscillators, where $x_{i} \in \mathbb R^{m}$ is the state of oscillator $i$ with limit cycle $\gamma_{i} \subset \real^{m}$ and period $T_{i}>0$. We assume that the oscillators are weakly coupled with interaction graph $G(\mc V,\mc E)$ and dynamics
\begin{equation}
	\dot x_{i} = f_{i}(x_{i}) + \epsilon \sum\nolimits_{\{i,j\} \in \mc E} g_{ij}(x_{i},x_{j})
	\,,\;\; i \in \until n \,,
	\label{eq: weakly coupled oscillator}
\end{equation}
where $g_{ij}(\cdot)$ is the coupling function for the pair $\{i,j\} \in \mc E$. The coupling $g_{ij}(\cdot)$ can possibly be impulsive. 
The weak coupling in \eqref{eq: weakly coupled oscillator} can be identified with the weak forcing in \eqref{eq: weakly forced oscillator}, and a transformation to phase coordinates yields 
\begin{equation*}
	\dot \varphi_{i} = \Omega_{i} + \epsilon \sum\nolimits_{\{i,j\} \in \mc E} Q_{i}(\varphi)g_{ij}(x_{i}(\varphi_{i}),x_{j}(\varphi_{j})) 
	\,,
\end{equation*}
where $\Omega_{i} = 2\pi/T_{i}$.
The local change of variables $\theta_{i}(t) = \varphi_{i}(t) - \Omega_{i} t$ then yields the coupled phase dynamics
\begin{equation*}
	\dot \theta_{i} =  \epsilon \sum\nolimits_{\{i,j\} \in \mc E} \!\! Q_{i}(\theta_{i}+\Omega_{i}t) g_{ij}(x_{i}(\theta_{i} + \Omega_{i}t),x_{j}(\theta_{j} + \Omega_{j}t)) 
	.
\end{equation*}
An averaging analysis applied to the $\theta$-dynamics results in
\begin{equation}
	\dot \theta_{i} =  \epsilon \omega_{i} + \epsilon \sum\nolimits_{\{i,j\} \in \mc E} h_{ij}(\theta_{i} - \theta_{j}) 
	\label{eq: averaged phase osc}
	\,,
\end{equation}
where $\omega_{i} = h_{ii}(0)$ and the averaged coupling functions are
\begin{equation*}
h_{ij}(\chi) = 
\lim_{T \to \infty} \frac{1}{T} \int_{0}^{T} \! Q_{i}(\Omega_{i}\tau) g_{ij}(x_{i}(\Omega_{i}\tau),x_{j}(\Omega_{j}\tau - \chi)) d\tau.
\end{equation*}
Notice that the averaged coupling functions $h_{ij}$ are $2\pi$-periodic and the coupling is diffusive. If all functions $h_{ij}$ are odd, a first-order Fourier series expansion of $h_{ij}$ yields $h_{ij}(\cdot) \approx a_{ij}\sin(\cdot)$ as first harmonic with some coefficient $a_{ij}$.  In this case, the  dynamics \eqref{eq: averaged phase osc} in the slow time scale $\tau = \epsilon t$ reduce exactly to the coupled oscillator model \eqref{eq: coupled oscillator model}.

This analysis justifies calling the coupled oscillator model \eqref{eq: coupled oscillator model} the {\em canonical model} for coupled limit-cycle oscillators.


\section{Synchronization Notions and Metrics}
\label{Section: Models and Synchronization Notions}

In this section, we introduce different notions of synchronization. Whereas the first four subsections address the commonly studied notions of synchronization associated with a coherent behavior and cohesive phases, Subsection \ref{Subsection: Phase Balancing} addresses the converse concept of phase balancing.

\subsection{Synchronization Notions}

The coupled oscillator model \eqref{eq: coupled oscillator model} evolves
on $\torus^{n}$, and features an important symmetry, namely the rotational
invariance of the angular variable $\theta$. This symmetry gives rise to
the rich synchronization dynamics. Different levels of synchronization can
be distinguished, and the most commonly studied notions are phase and
frequency synchronization.

{\bf Phase synchronization:} 
A solution $\theta:\, \real_{\geq 0} \to \torus^{n}$ to the coupled oscillator model \eqref{eq: coupled oscillator model} achieves {\em phase synchronization} if all phases $\theta_{i}(t)$ become identical as $t \to \infty$.

{\bf Phase cohesiveness:} 
As we will see later, phase synchronization can occur only if all natural frequencies $\omega_{i}$ are identical. If the natural frequencies are not identical, then each pairwise distance $|\theta_{i}(t) - \theta_{j}(t)|$ can converge to a constant but not necessarily zero value. The concept of phase cohesiveness formalizes this possibility. For $\gamma \in {[0,\pi[}$, let $\bar\Delta_{G}(\gamma) \subset \torus^{n}$ be the closed set of angle arrays $(\theta_{1},\dots,\theta_{n})$ with the property $|\theta_{i}-\theta_{j}| \leq \gamma$ for all $\{i,j\} \in \mc E$, that is, each pairwise phase distance is bounded by $\gamma$. Also, let $\Delta_{G}(\gamma)$ be the interior of $\bar\Delta_{G}(\gamma)$. Notice that $\bararc(\gamma) \subseteq \bar\Delta_{G}(\gamma)$ but the two sets are generally not equal. A solution $\theta:\, \real_{\geq 0} \to \torus^{n}$ is then said to be {\em phase cohesive} if there exists a length $\gamma \in {[0,\pi[}$ such that $\theta(t) \in \bar\Delta_{G}(\gamma)$ for all $t\geq 0$.

{\bf Frequency synchronization:}
 A solution $\theta:\, \real_{\geq 0} \to \torus^{n}$ achieves {\em frequency synchronization} if all frequencies $\dot\theta_{i}(t)$ converge to a common frequency $\subscr{\omega}{sync} \in \real$ as $t \to \infty$. The explicit synchronization frequency $\subscr{\omega}{sync} \in \real$ of the coupled oscillator model \eqref{eq: coupled oscillator model} can be obtained by summing over all equations in \eqref{eq: coupled oscillator model} as $\sum_{i=1}^{n} \dot \theta_{i} = \sum_{i=1}^{n} \omega_{i}$. In the frequency-synchronized case, this sum simplifies to $\sum_{i=1}^{n} \subscr{\omega}{sync} = \sum_{i=1}^{n} \omega_{i}$. In conclusion, if a solution of the coupled oscillator model \eqref{eq: coupled oscillator model} achieves frequency synchronization, then it does so with synchronization frequency equal to $\subscr{\omega}{sync} = \sum_{i=1}^{n} \omega_{i}/n$. By transforming to a rotating frame with frequency $\subscr{\omega}{sync}$~and by replacing $\omega_{i}$ by $\omega_{i} - \subscr{\omega}{sync}$, we obtain $\subscr{\omega}{sync} = 0$ (or~equivalently $\omega \in \fvec 1_{n}^{\perp}$). In what follows, without loss of generality, we will sometimes assume that $\omega \in \fvec 1_{n}^{\perp}$ so that $\subscr{\omega}{sync} = 0$.

\begin{remark}[\bf Terminology]
Alternative terminologies for phase synchronization include full, exact, or
perfect synchronization. For a frequency-synchronized solution all phase
distances $|\theta_{i}(t) - \theta_{j}(t)|$ are constant in a rotating
coordinate frame with frequency $\subscr{\omega}{sync}$, and the
terminology {\it phase locking} is sometimes used instead of frequency
synchronization. Other commonly used terms include frequency locking,
frequency entrainment, or also partial synchronization.  \oprocend
\end{remark}

{\bf Synchronization:} The main object under study in most applications and
theoretic analyses are phase cohesive and frequency-synchronized solutions,
that is, all oscillators rotate with the same synchronization frequency,
and all their pairwise phase distances are bounded. In the following, we
restrict our attention to synchronized solutions with sufficiently small
phase distances $|\theta_{i} - \theta_{j}| \leq \gamma < \pi/2$ for
$\{i,j\} \in \mc E$. Of course, there may exist other possible solutions,
but these are not necessarily stable (see our analysis in Section
\ref{Section: Analysis of Synchronization}) or not relevant in most
applications\footnote{For example, in power network applications the
  coupling terms $a_{ij} \sin(\theta_{i} - \theta_{j})$ are power flows
  along transmission lines $\{i,j\} \in \mc E$, and the phase distances
  $|\theta_{i} - \theta_{j}|$ are bounded well below $\pi/2$ due to thermal
  constraints. In Subsection \ref{Subsection: Phase Balancing}, we present a converse synchronization notion, where the goal is to maximize phase distances.}.
We say that a solution $\map{\theta}{\real_{\geq0}}{\torus^{n}}$ to the
coupled oscillator model \eqref{eq: coupled oscillator model} is {\em
  synchronized}
if there exists $\subscr{\theta}{sync} \in \bar\Delta_{G}(\gamma)$ for some
$\gamma \in {[0,\pi/2[}$ and $\subscr{\omega}{sync}\in\real$ (identically
    zero for $\omega \in \fvec 1_{n}^{\perp}$) such that $\theta(t) =
    \subscr{\theta}{sync} + \subscr{\omega}{sync}\fvec 1_n t \pmod{2\pi}$
    for all $t \geq 0$.

\begin{figure}[htbp]
	\centering{
	\includegraphics[width = 0.85\columnwidth]{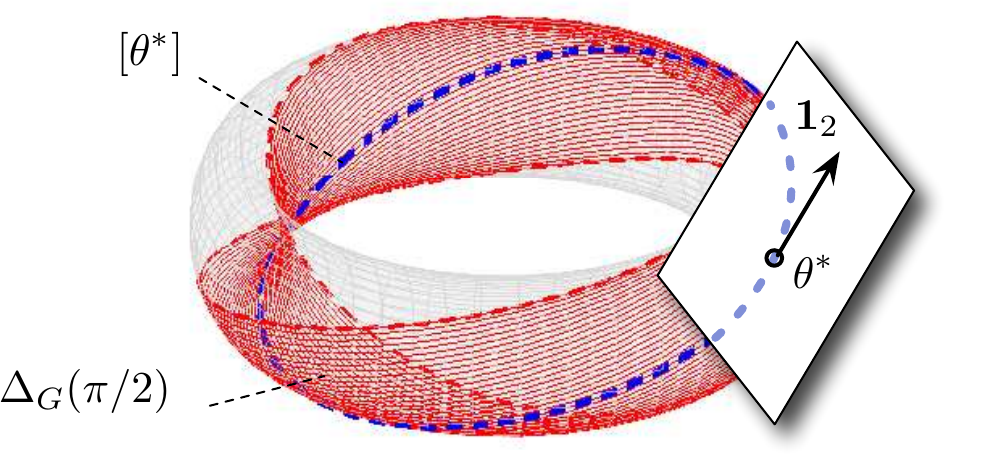}
	\caption{Illustration of the state space $\torus^{2}$, the set $\Delta_{G}(\pi/2)$, the synchronization manifold $[\theta^{*}]$ associated to a phase-synchronized angle array $\theta^{*} =$ $ (\theta_{1}^{*},\theta_{2}^{*}) \in \bar\Delta_{G}(0)$, and the tangent space with translation vector $\fvec 1_{2}$ at~$\theta^{*}$.}
	\label{Fig: sync manifold}
	}
\end{figure}

{\bf Synchronization manifold:} The geometric object under study in
synchronization is the synchronization manifold.
Given a point $r\in\mycircle$ and an angle $s\in[0,2\pi]$, let
$\rot_s(r)\in\mycircle$ be the rotation of $r$ counterclockwise by the
angle $s$.  For $(r_{1},\dots,r_{n}) \in \torus^n $, define the equivalence
class
\begin{equation*}
  [(r_1,\dots,r_n)] \!=\! \setdef{ ( \rot_s(r_1), \dots, \rot_s(r_n) ) \in \mbb
    T^n\! }{ \! s\in[0,2\pi]}.
\end{equation*}
Clearly, if $(r_1,\dots,r_n)\in\bar\Delta_G(\gamma)$ for some $\gamma \in
{[0,\pi/2[}$, then $[(r_1,\dots,r_n)]\subset\bar\Delta_G(\gamma)$.
Given a synchronized solution characterized by
$\subscr{\theta}{sync}\in\bar\Delta_G(\gamma)$ for some $\gamma \in
{[0,\pi/2[}$, the set $[\subscr{\theta}{sync}] \subset
    \bar\Delta_G(\gamma)$ is a {\em synchronization manifold} of the
    coupled-oscillator model \eqref{eq: coupled oscillator model}.
Note that a synchronized solution takes value in a synchronization manifold
due to rotational symmetry, and for $\omega \in \fvec 1_{n}^{\perp}$
(implying $\subscr{\omega}{sync} = 0$) a synchronization manifold is also
an equilibrium manifold of the coupled oscillator model \eqref{eq: coupled
  oscillator model}.
These geometric concepts are illustrated in Figure \ref{Fig: sync manifold}
for the two-dimensional case.


\subsection{A Simple yet Illustrative Example}
\label{Subsection: 2d example}

\begin{figure}[h]
    \centering{
    \subfigure[Vector field \eqref{eq: difference dynamics for two oscillators} for  $\theta_{2} - \theta_{1} > 0$]
    {
    \includegraphics[width = 0.45\columnwidth]{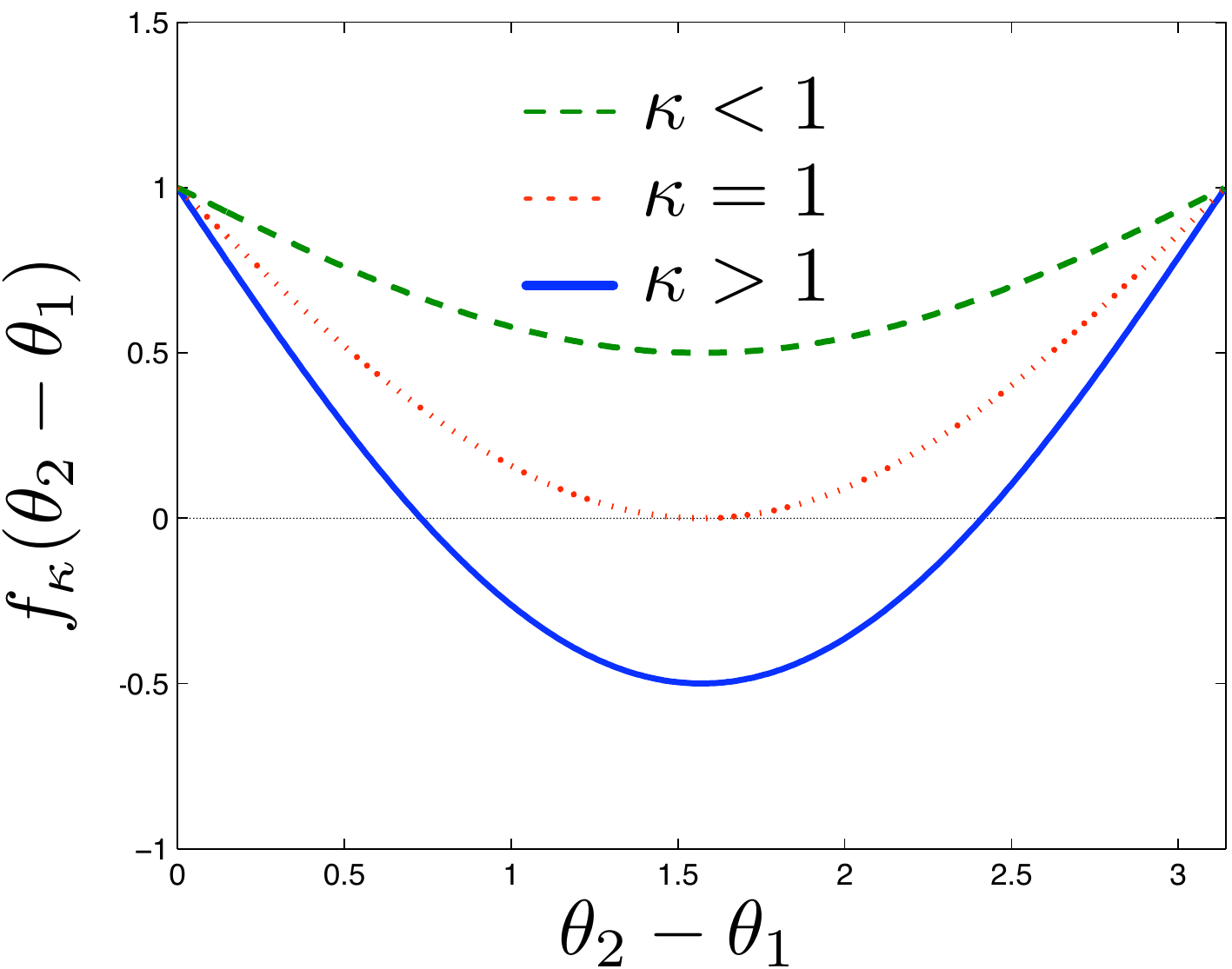}
    \label{Fig: vector field for two oscillators}
    }
   \hspace{0.15cm}
    \subfigure[Trajectory $\theta(t)$ for $\kappa = 1$]
    {
    \includegraphics[width = 0.45\columnwidth]{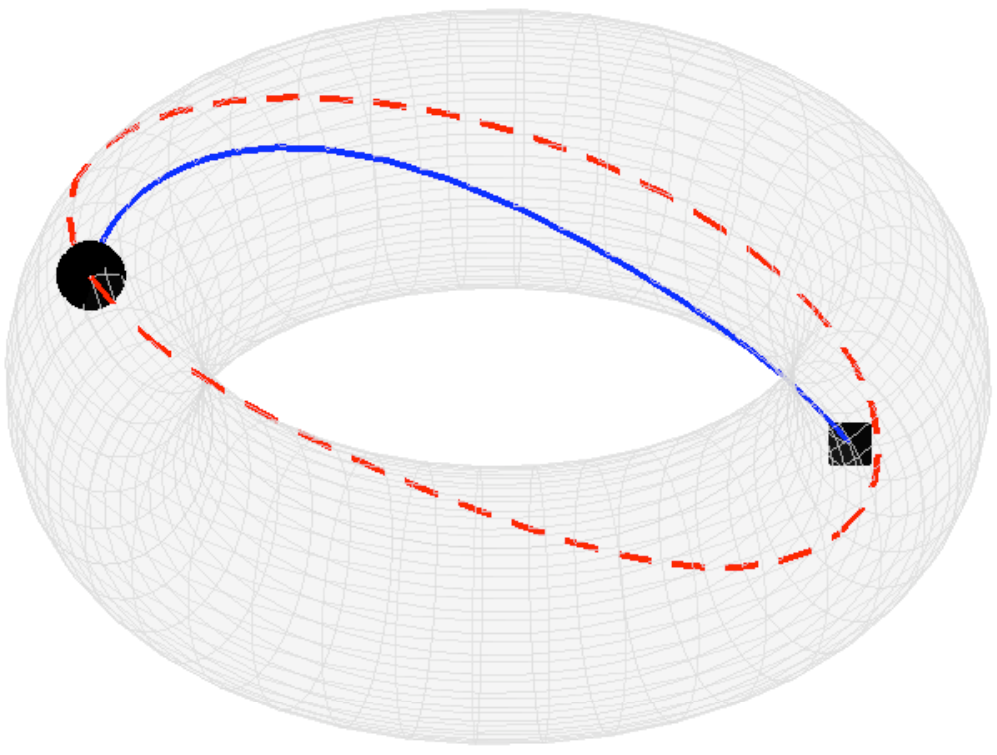}
    \label{Fig: evolution in configuration space}
    }
    \caption{Plot of the vector field \eqref{eq: difference dynamics for
        two oscillators} for various values of $\kappa$ and a trajectory
      $\theta(t) \in \torus^{2}$ for the critical case $\kappa = 1$,
      where the dashed line is the saddle equilibrium manifold and
      $\blacksquare$ and {\large$\bf\bullet$} depict $\theta(0)$ and
      $\lim_{t \to \infty}\theta(t)$. The non-smoothness of the vector
      field $f(\theta_{2} - \theta_{1})$ at the boundaries $\{0,\pi\}$ is
      an artifact of the non-smoothness of the geodesic distance on $\mbb
      T^{2}$}
    \label{Fig: two oscillators}
    }
\end{figure}

The following example illustrates the different notions of synchronization
introduced above and points out various important geometric subtleties
occurring on the compact state space $\torus^{2}$.
Consider $n=2$ oscillators with $\omega_{2} \geq 0 \geq \omega_{1} =
-\omega_{2}$. We restrict our attention to angles contained in an open
half-circle: for angles $\theta_{1}$, $\theta_{2}$ with $|\theta_{2} -
\theta_{1}|<\pi$, the {\it angular difference} $\theta_{2}-\theta_{1}$ is
the number in ${]\!-\!\pi,\pi[}$ with magnitude equal to the geodesic
distance $|\theta_{2} - \theta_{1}|$ and with positive sign if and only if
the counter-clockwise path length from $\theta_{1}$ to $\theta_{2}$ on
$\mbb T^{1}$ is smaller than the clockwise path length. With this
definition the two-dimensional oscillator dynamics $(\dot \theta_{1},\dot
\theta_{2})$ can be reduced to the scalar difference dynamics $\dot
\theta_{2} - \dot \theta_{1}$. After scaling time as $t \mapsto t
(\omega_{2} - \omega_{1})$ and introducing $\kappa = 2 a_{12}/(\omega_{2} -
\omega_{1})$ the difference dynamics are
\begin{equation}
	\dt (\theta_{2} - \theta_{1} )
	=
	f_{\kappa}(\theta_{2} - \theta_{1})
	:=
	1 - \kappa \sin(\theta_{2} - \theta_{1})
	\label{eq: difference dynamics for two oscillators}
	\,.
\end{equation}
The scalar dynamics \eqref{eq: difference dynamics for two oscillators} can be analyzed graphically by plotting the vector field $f_{\kappa}(\theta_{2} - \theta_{1})$ over the difference variable $\theta_{2} - \theta_{1}$, as in Figure \ref{Fig: vector field for two oscillators}. 
Figure \ref{Fig: vector field for two oscillators} displays a saddle-node bifurcation at $\kappa=1$. For $\kappa < 1$ no equilibrium of \eqref{eq: difference dynamics for two oscillators} exists, and for $\kappa >1$ an asymptotically stable equilibrium $\subscr{\theta}{stable} = \arcsin(\kappa^{-1}) \in {]0,\pi/2[}$ together with a saddle point $\subscr{\theta}{saddle} = \arcsin(\kappa^{-1}) \in {]\pi/2,\pi[}$ exists. 

For $\theta(0) \in \arc(|\subscr{\theta}{saddle}|)$ all trajectories converge exponentially to $\subscr{\theta}{stable}$, that is, the oscillators synchronize exponentially. Additionally, the oscillators are phase cohesive if an only if $\theta(0) \in \bararc(|\subscr{\theta}{saddle}|)$, where all trajectories remain bounded. For $\theta(0) \not \in \bararc(|\subscr{\theta}{saddle}|)$ the difference $\theta_{2}(t) - \theta_{1}(t)$ will increase beyond $\pi$, and by definition will change its sign since the oscillators change orientation. Ultimately, $\theta_{2}(t) - \theta_{1}(t)$ converges to the equilibrium $\subscr{\theta}{stable}$ in the branch where $\theta_{2} - \theta_{1} < 0$. In the configuration space $\torus^{2}$ this implies that the distance $\abs{\theta_{2}(t) - \theta_{1}(t)}$ increases to its maximum value $\pi$ and shrinks again, that is, the oscillators are not phase cohesive and revolve once around the circle before converging to the equilibrium manifold. Since $\sin(\subscr{\theta}{stable}) = \sin(\subscr{\theta}{saddle}) = \kappa^{-1}$, strongly coupled oscillators with $\kappa \gg 1$ practically achieve phase synchronization from every initial condition in an open semi-circle. 
In the critical case, $\kappa = 1$, the saddle equilibrium manifold at
$\pi/2$ is globally attractive but not stable. An representative trajectory
is illustrated in Figure~\ref{Fig: evolution in configuration space}.

In conclusion, the simple but already rich $2$-dimensional case shows that
two oscillators are phase cohesive and synchronize if and only if $\kappa >
1$, that is, if and only if the coupling dominates the non-uniformity as
$2a_{12} > \omega_{2} - \omega_{1}$. The ratio $1/\kappa$ determines the
ultimate phase cohesiveness as well as the set of admissible initial
conditions. For $\kappa \gg 1$, practical phase synchronization is achieved
for all angles in an open semi-circle.
More general coupled oscillator networks display the same phenomenology,
but the threshold from incoherence to synchrony is generally unknown.

\subsection{Synchronization Metrics}

The notion of phase cohesiveness can be understood as a performance measure
for synchronization and phase synchronization is simply the extreme case of
phase cohesiveness with $\lim_{t \to \infty} \theta(t) \in
\bar\Delta_{G}(0) = \bararc(0)$.
An alternative performance measure is the magnitude of the so-called {\em
  order parameter} introduced by Kuramoto \cite{YK:75,YK:84}:
\begin{equation*}
r e^{\mathrm{i} \psi} = \frac{1}{n} \sum\nolimits_{j=1}^{n} e^{\mathrm{i} \theta_{j}}
\,.
\end{equation*}
The order parameter
is the centroid of all oscillators represented as points on the unit circle
in $\mbb C^{1}$. The magnitude $r$ of the order parameter is a
synchronization measure: if all oscillators are phase-synchronized, then
$r=1$, and if all oscillators are spaced equally on the unit circle, then $r=0$.
The latter case is characterized in Subsection~\ref{Subsection: Phase Balancing}.~
For~a~complete graph, the magnitude $r$ of the order parameter serves as an
{\it average} performance index for synchronization, and phase cohesiveness
can be understood as a\,{\it worst-case} performance index. Extensions of
the order parameter tailored to non-complete graphs have been proposed
in \cite{AJ-NM-MB:04,DAP-NEL-RS-DG-JKP:07,LS-AS-RS:06}.

\begin{figure}[htbp]
  \centering{
    \includegraphics[scale=0.4]{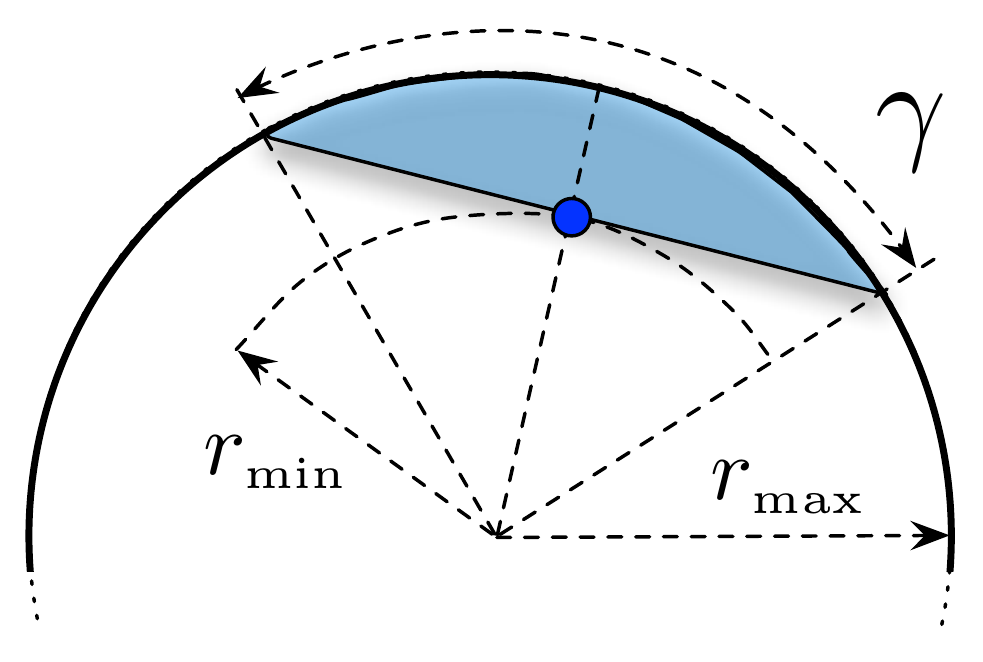}
    \caption{Schematic illustration of an arc of length $\gamma \in
      {[0,\pi]}$, its convex hull (shaded), and the value
      {\large\bf${\color{blue}\fvec\bullet}$} of the corresponding order
      parameter $r e^{\mathrm{i}\psi}$ with minimum magnitude
      $\subscr{r}{min} =\cos(\gamma/2)$ and maximum magnitude $\subscr{r}{max} = 1$.}
    \label{Fig: arc and its convex hull}
  }
\end{figure}

For a complete graph and for $\gamma$ sufficiently small, the set
$\bar\Delta_{G}(\gamma)$ reduces to $\bararc(\gamma)$, the arc of length
$\gamma$ containing all oscillators. The order parameter is contained
within the convex hull of this arc since it is the centroid of all
oscillators, see Figure \ref{Fig: arc and its convex hull}. In this case,
the magnitude $r$ of the order parameter can be related to the arc length
$\gamma$.

\begin{lemma}{(\bf Shortest arc length and order parameter, \cite[Lemma 2.1]{FD-FB:10w})} 
\label{Lemma: phase cohesiveness and order parameter}
  Given an angle array $\theta = (\theta_{1},\dots,\theta_{n}) \in \torus^n$ with $n
  \geq 2$, let $r(\theta) = \frac{1}{n} |\sum_{j=1}^{n}
  e^{\mathrm{i}\theta_{j}} |$ be the magnitude of the order parameter, and
  let $\gamma(\theta)$ be the length of the shortest arc containing all
  angles, that is, $\theta\in\bararc(\gamma(\theta))$.  The following
  statements hold:
  \begin{enumerate}
  \item[1)] if $\gamma(\theta) \in {[0,\pi]}$, then $r(\theta) \in
    [\cos(\gamma(\theta)/2), 1]$; and

  \item[2)] if $\theta \in\bararc(\pi)$, then $\gamma(\theta) \in
    [2\arccos(r(\theta)), \pi]$.
\end{enumerate}
\end{lemma}

\subsection{Synchronization Conditions}
\label{Subsection: Synchronization Conditions}

The coupled oscillator dynamics \eqref{eq: coupled oscillator model}
feature (i) the synchronizing coupling described by the graph $G(\mc V,\mc
E,A)$ and (ii) the de-synchronizing effect of the non-uniform natural
frequencies $\omega$.  Loosely speaking, synchronization occurs when the
coupling dominates the non-uniformity.  Various conditions have been
proposed in the synchronization and power systems literature to quantify
this trade-off.

The coupling is typically quantified by the algebraic connectivity
$\lambda_{2}(L)$
\cite{FFW-SK:80,FD-FB:09z,AJ-NM-MB:04,TN-AEM-YCL-FCH:03,AA-ADG-JK-YM-CZ:08,SB-VL-YM-MC-DUH:06}
or the weighted nodal degree $\textup{deg}_{i} \triangleq
\sum\nolimits_{j=1}^{n} a_{ij}$
\cite{FW-SK:82,FD-FB:11d,GK-MBH-KEB-MJB-BK-DA:06,FD-FB:09z,LB-LS-ADG:09},
and the non-uniformity is quantified by either absolute norms $\| \omega
\|_{p}$ or incremental norms $\| B^{T} \omega \|_{p}$, where typically $p
\in \{2,\infty\}$. Sometimes, these conditions can be evaluated only
numerically since they are state-dependent \cite{FFW-SK:80,FW-SK:82} or
arise from a non-trivial linearization process, such as the Master
stability function formalism
\cite{AA-ADG-JK-YM-CZ:08,SB-VL-YM-MC-DUH:06,LMP-TLC:98}. In general,
concise and accurate results are known\,only for specific topologies such as
 complete graphs \cite{FD-FB:10w}, linear chains \cite{SHS-REM:88}, and
bipartite graphs \cite{MV-OM:09} with uniform weights.

For arbitrary coupling topologies only sufficient conditions are known
\cite{FFW-SK:80,FD-FB:09z,AJ-NM-MB:04,FW-SK:82} as well as numerical
investigations for random networks
\cite{JGG-YM-AA:07,JCB-LDV-MJP:11,TN-AEM-YCL-FCH:03,YM-AFP:04,ACK:10}. Simulation
studies indicate that these conditions are conservative estimates on the
threshold from incoherence to synchrony. Literally, every review article on
synchronization draws attention to the problem of finding sharp
synchronization conditions
\cite{SHS:01,JAA-LLB-CJPV-FR-RS:05,FD-FB:10w,SHS:00,AA-ADG-JK-YM-CZ:08,SB-VL-YM-MC-DUH:06,FD-MC-FB:11v-arxiv}.

\subsection{Phase Balancing and Splay State}
\label{Subsection: Phase Balancing}

In certain applications in neuroscience \cite{FV-JPL-ER-JM:01,EB-PH-JM:03,SMC-GBE-MCV-JMB:97}, deep-brain stimulation \cite{PAT:03,AN-JM:11}, and vehicle coordination \cite{DAP-NEL-RS-DG-JKP:07,RS-DP-NEL:07,RS-DP-NEL:08,DJK:08,DJK-PL-KAM-TJ:08}, one is not interested in the coherent behavior with synchronized (or nearly synchronized) phases, but rather in the phenomenon of synchronized frequencies and de-sychronized phases. 

Whereas the phase-synchronized state is characterized by the order parameter $r$ achieving its maximal (unit) magnitude, we say that a solution $\theta:\, \real_{\geq 0} \to \torus^{n}$ to the coupled oscillator model \eqref{eq: coupled oscillator model} achieves {\em phase balancing} if all phases $\theta_{i}(t)$ converge to $\bal = \{\theta \in \torus^{n}:\, r(\theta) = |\frac{1}{n} \sum\nolimits_{j=1}^{n} e^{\mathrm{i} \theta_{j}}| = 0\}$ as $t \to \infty$,
that is, the oscillators are distributed over the unit circle $\mbb S^{1}$, such that their centroid $r e^{\mathrm{i} \psi}$ vanishes. 
We refer to \cite{RS-DP-NEL:07} for a geometric characterization of the balanced~state.

One balanced state of particular interest in neuroscience applications \cite{FV-JPL-ER-JM:01,EB-PH-JM:03,SMC-GBE-MCV-JMB:97,PAT:03,AN-JM:11} is the so-called splay state $\{\theta \in \torus^{n}:\, \theta_{i} = {i \cdot 2\pi/n + \varphi}\pmod{2\pi} \,, \varphi \in \mbb S^{1} \,, i \in \until n  \} \subseteq \bal$ corresponding to phases uniformly distributed around the unit circle $\mbb S^{1}$ with distances $2\pi/n$. Other highly symmetric balanced states consist of multiple clusters of collocated phases, where the clusters themselves are arranged in splay state, see \cite{RS-DP-NEL:07,RS-DP-NEL:08}.


\section{Analysis of Synchronization}
\label{Section: Analysis of Synchronization}

In this section we present several analysis approaches to synchronization
in the coupled oscillator model \eqref{eq: coupled oscillator model}. We
begin with a few basic ideas to provide important intuition as well as the
analytic basis for further analysis.

\subsection{Some Simple Yet Important Insights}
\label{Subsection: simple insights}

The potential energy $\map{U}{\torus^{n}}{\real}$ of the elastic spring
network in Figure \ref{Fig: Mechanical analog} is, up to an additive
constant, given~by
\begin{equation}
	U(\theta) = \sum\nolimits_{\{i,j\} \in \mc E} a_{ij} \bigl( 1 - \cos(\theta_{i} - \theta_{j}) \bigr)
	\,.
	\label{eq: potential energy}
\end{equation}
By means of the potential energy, the coupled oscillator model \eqref{eq:
  coupled oscillator model} can reformulated as the forced gradient system
\begin{equation}
	\dot \theta_{i}
	=
	\omega_{i} - \nabla_{i} U(\theta)
	\,, \qquad i \in \until n
	\label{eq: coupled oscillator model - potential}
	\,,
\end{equation}
where $\nabla_{i} U(\theta) = \frac{\partial}{\partial \theta_{i}}
U(\theta)$ denotes the partial derivative. It can be easily verified that
the phase-synchronized state $\theta_{i} = \theta_{j}$ for all $\{i,j\} \in
\mc E$ is a local minimum of the potential energy \eqref{eq: potential
  energy}. The gradient formulation \eqref{eq: coupled oscillator model -
  potential} clearly emphasizes the competition between the
synchronization-enforcing coupling through the potential $U(\theta)$ and
the synchronization-inhibiting heterogeneous natural frequencies
$\omega_{i}$.

We next note that $\omega$ has to satisfy certain bounds, relative to the
weighted nodal degree, in order for a synchronized solution to exist.

\begin{lemma}{\bf(Necessary sync conditions)}
\label{Lemma: Necessary sync condition}
Consider the coupled oscillator model \eqref{eq: coupled oscillator model}
with graph $G(\mc V,\mc E,A)$, frequencies $\omega \in \fvec
1_{n}^{\perp}$, and nodal degree $\textup{deg}_{i} =
\sum\nolimits_{j=1}^{n} a_{ij}$ for each oscillator $i \in \until n$. If there
exists a synchronized solution $\theta \in \bar\Delta_{G}(\gamma)$ for some
$\gamma \in {[0,\pi/2]}$, then the following conditions hold:
\begin{enumerate}
\item[1)] {\bf Absolute bound:} 
  For each node $i \in \until n$,
  \begin{equation}
    \textup{deg}_{i} \sin(\gamma) 
    \geq |\omega_{i}| \,;
    \label{eq: necessary sync condition - 1}
  \end{equation} 

\item[2)] {\bf Incremental bound:} 
  For all distinct $i,j \in \until n$,
  \begin{equation}
    (\textup{deg}_{i} + \textup{deg}_{j})\sin(\gamma) 
    \geq |\omega_{i} - \omega_{j}| \,. 
    \label{eq: necessary sync condition - 2}
  \end{equation} 
\end{enumerate}
\end{lemma}

\smallskip\begin{proof}
Since $\omega \in \fvec 1_{n}^{\perp}$, the synchronization frequency
$\subscr{\omega}{sync}$ is zero, and phase and frequency synchronized
solutions are equilibrium solutions determined by the equations
\begin{equation}
	\omega_{i} = \sum\nolimits_{j=1}^{n} a_{ij} \sin(\theta_{i} - \theta_{j})
	\,, \qquad i \in \until n \,.
	\label{eq: equilibrium equations - component-wise}
\end{equation}
Since $\sin(\theta_{i} - \theta_{j}) \in {[-\sin(\gamma),+\sin(\gamma)]}$ for $\theta \in \bar\Delta_{G}(\gamma)$, the equilibrium equations \eqref{eq: equilibrium equations - component-wise} have no solution if condition \eqref{eq: necessary sync condition - 1} is not satisfied. Since $\omega \in \fvec 1_{n}^{\perp}$, an incremental bound on $\omega$ seems to be more appropriate than an absolute bound. The  subtraction of the $i$th and $j$th equation \eqref{eq: equilibrium equations - component-wise} yields
\begin{equation*}
	\omega_{i} - \omega_{j}
	=
	\sum\nolimits_{k=1}^{n} \left( a_{ik} \sin(\theta_{i}-\theta_{k}) - a_{jk} \sin(\theta_{j}-\theta_{k}) \right)
	\,.
	\label{eq: eq: fixed-point equations - component-wise - difference}
\end{equation*}
Again, since the coupling is bounded, the above equation has no solution in $\bar\Delta_{G}(\gamma)$ if condition \eqref{eq: necessary sync condition - 2} is not satisfied.
\end{proof}

The following result is fundamental for various approaches to phase and
frequency synchronization. To the best of the authors' knowledge this result
has been first established in \cite{CJT-OJMS:72a}, and it has been reproved 
numerous times.

\begin{lemma}{(\bf Stable synchronization in $\Delta_{G}(\pi/2)$)}
\label{Lemma: stable equilibrium in pi/2 arc}
Consider the coupled oscillator model \eqref{eq: coupled oscillator model}
with a connected graph $G(\mc V,\mc E,A)$ and frequencies $\omega \in \fvec
1_{n}^{\perp}$. The following statements hold:
\begin{enumerate}
\item[1)] {\bf Jacobian:} The Jacobian $J(\theta)$ of the coupled
  oscillator model \eqref{eq: coupled oscillator model} evaluated at
  $\theta \in \torus^{n}$ is given by
  \begin{equation*}
    J (\theta) 
    =
    - B \diag(\{a_{ij}\cos(\theta_{i}-\theta_{j})\}_{\{i,j\} \in \mc E}) B^{T}
    \,;
  \end{equation*}

\item[2)] {\bf Local stability and uniqueness:} If there exists an
  equilibrium $\theta^{*} \in \Delta_{G}(\pi/2)$, then
      \begin{enumerate}
        
      \item[(i)] $-J(\theta^{*})$ is a  Laplacian matrix;
        
      \item[(ii)] the equilibrium manifold $[\theta^{*}] \in
        \Delta_{G}(\pi/2)$ is locally exponentially stable; and
	
      \item[(iii)] this equilibrium manifold is unique in
        $\bar\Delta_{G}(\pi/2)$.
      \end{enumerate}
      
\end{enumerate}
\end{lemma}

\smallskip\begin{proof} Since $\frac{\partial}{\partial \theta_{i}} \bigl(
\omega_{i} - \sum\nolimits_{j=1}^{n} a_{ij} \sin(\theta_{i}-\theta_{j})
\bigr) = -\sum\nolimits_{j=1}^{n} a_{ij} \cos(\theta_{i} - \theta_{j})$ and
$\frac{\partial}{\partial \theta_{j}} \bigl( \omega_{i} -
\sum\nolimits_{j=1}^{n} a_{ij} \sin(\theta_{i}-\theta_{j}) \bigr) =
a_{ij}\cos(\theta_{i} - \theta_{j})$, we obtain that the Jacobian is equal
to minus the Laplacian matrix of the connected graph $G(\mc V,\mc E,\tilde
A)$ with the (possibly negative) weights $\tilde a_{ij} =
a_{ij}\cos(\theta_{i}-\theta_{j})$. Equivalently, in compact notation $J
(\theta) = - B \diag(\{a_{ij}\cos(\theta_{i}-\theta_{j})\}_{\{i,j\} \in \mc
  E}) B^{T}$. This completes the proof of statement 1).

The Jacobian $J(\theta)$ evaluated for an equilibrium $\theta^{*} \in
\Delta_{G}(\pi/2)$ is minus the Laplacian matrix of the graph $G(\mc V,\mc
E,\tilde A)$ with strictly positive weights $\tilde a_{ij} =
a_{ij}\cos(\theta_{i}^{*}-\theta_{j}^{*}) > 0$ for every $\{i,j\} \in \mc
E$. Hence, $J(\theta^{*})$ is negative semidefinite with the nullspace
$\fvec 1_{n}$ arising from the rotational symmetry, see Figure \ref{Fig:
  sync manifold}.
Consequently, the equilibrium point $ \theta^{*} \in \Delta_{G}(\pi/2)$ is
locally (transversally) exponentially stable, or equivalently, the
corresponding equilibrium manifold $[\theta^{*}] \in \Delta_{G}(\pi/2)$ is
locally exponentially stable.

The uniqueness statement follows since the right-hand side of the coupled
oscillator model \eqref{eq: coupled oscillator model} is a one-to-one
function (modulo rotational symmetry) for $\theta \in \bar\Delta_{G}(\pi/2)$, see \cite[Corollary
  1]{AA-SS-VP:81}.
This completes the proof of statement 2).
\end{proof}

By Lemma \ref{Lemma: stable equilibrium in pi/2 arc}, any equilibrium in $\Delta_{G}(\pi/2)$ is stable which supports the notion of phase cohesiveness as a performance metric.
Since the Jacobian $J(\theta)$ is the negative Hessian of the potential $U(\theta)$ defined in \eqref{eq: potential energy}, Lemma \ref{Lemma: stable equilibrium in pi/2 arc} also implies that any equilibrium in $\Delta_{G}(\pi/2)$ is a local minimizer of $U(\theta)$.
Of particular interest are so-called {\em $\sphere^{1}$-synchronizing
  graphs} for which all critical points of \eqref{eq: potential energy} are
hyperbolic, the phase-synchronized state is the global minimum of
$U(\theta)$, and all other critical points are local maxima or saddle
points. The class of $\sphere^{1}$-synchronizing graphs includes, among
others, complete graphs and acyclic
graphs~\cite{AS:09,EAC-PAM-RF:10-a,EAC-PAM-RF:10-b,PM:06}.
 
These basic insights motivated various characterizations and explorations
of the critical points and the curvature of the potential $U(\theta)$ in
the literature on synchronization
\cite{AS:09,EAC-PAM-RF:10-a,EAC-PAM-RF:10-b,PM:06,PM:06,EM-AT:10,AJ-NM-MB:04,FD-FB:09z,FD-FB:10w,JCB-LDV-MJP:11,AS:09}
as well as on power systems
\cite{CJT-OJMS:72a,CJT-OJMS:72b,SS-PV:80,ARB-DJH:81,MI:92,AA-SS-VP:81,FW-SK:82,FFW-SK:80,KSC-DJH:86}.

\subsection{Phase Synchronization}

If all natural frequencies are identical, $\omega_{i} \equiv \omega$ for all $i \in \until n$, then a transformation of the coupled oscillator model \eqref{eq: coupled oscillator model} to a rotating frame with frequency $\omega$ leads to
\begin{equation}
	\dot \theta_{i}
	=
	- \sum\nolimits_{j=1}^{n} a_{ij} \sin(\theta_{i}-\theta_{j})
	\,, \qquad i \in \until n
	\label{eq: coupled oscillator model - phase sync}
	\,.
\end{equation}
The analysis of the coupled oscillator model \eqref{eq: coupled oscillator model - phase sync} is particularly simple and local phase synchronization can be concluded  by various analysis methods. A sample of different analysis schemes (by far not complete) includes the contraction property \cite{ZL-BF-MM:07,AS:09,UM-AP-FA:09,FD-FB:09z,RS-AS-PR:10}, quadratic\,Lyapunov functions \cite{AJ-NM-MB:04,FD-FB:09z}, linearization \cite{EC-PM:08,PM:06}, or order parameter and potential function arguments \cite{LS-AS-RS:06,RS-DP-NEL:07,JL:11}. 

The following theorem on phase synchronization summarizes a collection of
results originally presented in
\cite{LS-AS-RS:06,ZL-BF-MM:07,PM:06,AS:09,RS-DP-NEL:07,FD-FB:10w}, and it
can be easily proved given the insights developed in Subsection
\ref{Subsection: simple insights}.


\begin{theorem}{\bf(Phase synchronization)}
\label{Theorem: Phase sync}%
Consider the coupled oscillator model \eqref{eq: coupled oscillator model}
with a connected graph $G(\mc V,\mc E,A)$ and with frequency $\omega \in
\real^{n}$ (not necessarily zero mean).
The following statements are equivalent:
\begin{enumerate}
\item[(i)]{\bf Stable phase sync:} there exists a locally exponentially
  stable phase-synchronized solution $\theta \in \bararc(0)$ (or a synchronization manifold $[\theta] \in
  \bar\Delta_{G}(0)$); and
	
\item[(ii)]{\bf Uniformity:} there exists a constant $\omega \in \real$
  such that $\omega_{i} = \omega$ for all $i \in \until n$.
\end{enumerate}
If the two equivalent cases (i) and (ii) are true, the following statements hold:
\begin{enumerate}
	
\item[1)] {\bf Global convergence:} For all initial angles $\theta(0) \in
  \torus^{n}$ all frequencies $\dot\theta_{i}(t)$ converge to $\omega$ and
  all phases $\theta_{i}(t) - \omega t \pmod{2\pi}$ converge to the
  critical points $\{\theta \in \torus^{n}:\, \nabla U(\theta) = \fvec
  0_{n} \}$;
	
\item[2)] {\bf Semi-global stability:} The region of attraction of the
  phase-synchronized solution $\theta \in \bararc(0)$ contains the open
  semi-circle $\arc(\pi)$, and each arc $\bararc(\gamma)$ is
  positively invariant for every arc length $\gamma < \pi$;
  
  \item[3)] {\bf Explicit phase:} For initial angles in an open semi-circle
    $\theta(0) \in \arc(\pi)$, the asymptotic\,synchronization phase\,is\,%
    given\,by%
  \footnote{This ``average'' of angles (points on $\mathbb S^{1}$) is
    well-defined in an open semi-circle. If the parametrization of $\theta$
    has no discontinuity inside the arc containing all angles, then the
    average can be obtained by the usual formula.  }  $\theta(t) \!=\!
  \sum_{i=1}^{n}\! \theta_{i}(0)/n + \omega t \pmod{2\pi}$;
		
\item[4)] {\bf Convergence rate:} For every initial angle $\theta(0) \in
  \bararc(\gamma)$ with $\gamma < \pi$, the exponential convergence rate
  to phase synchronization is no worse than $\subscr{\lambda}{ps} = -
  \lambda_{2}(L) \sinc(\gamma)$; and
	
\item[5)] {\bf Almost global stability:} If the graph $G(\mc V,\mc E,A)$ is
  $\sphere^{1}$-synchronizing, the region of attraction of the
  phase-synchronized solution $\theta \in \bararc(0)$ is almost all of
  $\torus^{n}$.
\end{enumerate} 
\end{theorem}

\smallskip\begin{proof} {\em Implication (i) $\!\implies\!$ (ii):} By
assumption, there exist constants $\subscr{\theta}{sync} \in \sphere^{1}$
and $\subscr{\omega}{sync} \in \real$ such that $\theta_{i}(t) =
\subscr{\theta}{sync} + \subscr{\omega}{sync}t \pmod{2\pi}$. In the
phase-synchronized case, the dynamics \eqref{eq: coupled oscillator model}
then read as $\subscr{\omega}{sync} = \omega_{i}$ for all $i \in \until
n$. Hence, a necessary condition for the existence of phase synchronization
is that all $\omega_{i}$ are identical.

{\em Implication (ii) $\implies$ (i):} Consider the model \eqref{eq:
  coupled oscillator model} written in a rotating frame with frequency
$\omega$ as in \eqref{eq: coupled oscillator model - phase sync}. Note that
the set of phase-synchronized solutions $\bar\Delta_{G}(0)$ is an
equilibrium manifold. By Lemma \ref{Lemma: stable equilibrium in pi/2 arc},
we conclude that $\bar\Delta_{G}(0)$ is locally exponentially stable. This
concludes the proof of (i) $\Leftrightarrow$ (ii).

{\em Statement 1):} Note that \eqref{eq: coupled oscillator model - phase
  sync} can be written as the gradient flow $\dot \theta = - \nabla
U(\theta)$, and the corresponding potential function $U(\theta)$ is
non-increasing along trajectories. Since the sublevel sets of $U(\theta)$
are compact and the vector field $\nabla U(\theta)$ is smooth, the
invariance principle \cite[Theorem 4.4]{HKK:02} asserts that every solution
converges to set of equilibria of \eqref{eq: coupled oscillator model -
  phase sync}.

{\em Statements 2):} The coupled oscillator model \eqref{eq: coupled
  oscillator model - phase sync} can be re-written as the consensus-type
system
\begin{equation}
  \dot \theta_{i} = - \sum\nolimits_{j=1}^{n} b_{ij}(\theta) \cdot (\theta_{i} - \theta_{j})
  \,,\qquad i \in \until n\,,
  \label{eq: coupled oscillator model - phase sync - 2}
\end{equation}
where the weights $b_{ij}(\theta) = a_{ij} \sinc(\theta_{i} - \theta_{j})$
depend explicitly on the system state. Notice that for $\theta \in
\bararc(\gamma)$ and $\gamma < \pi$ the weights $b_{ij}(\theta)$ are
upper and lower bounded as $b_{ij}(\theta) \in {[ a_{ij}
    \sinc(\gamma),a_{ij} ]}$
Assume that the initial angles $\theta_i(0)$ belong to the set
$\bararc(\gamma)$, that is, they are all contained in an arc of length
$\gamma \in {[0,\pi[}$. In this case, a natural Lyapunov function to
    establish phase synchronization can be obtained from the {\em
      contraction property}, which aims at showing that the convex hull
    containing all oscillators is decreasing, see
    \cite{LM:08-arxiv,ZL-BF-MM:07,AS:09,UM-AP-FA:09,FD-FB:09z} and
    the\,review \cite[Section 2]{RS-AS-PR:10}.

Recall the geodesic distance between two angles on $\mathbb{S}^1$ and
define the continuous function
$\map{V}{\mathbb{T}^n}{[0,\pi]}$ by
\begin{equation} 
  V(\psi) =  \max\setdef{|\psi_i-\psi_j|}{i,j\in\until{n}}.
  \label{eq: contraction Lyapunov function}
\end{equation}
Notice that, if all angles are contained in an arc at time $t$, then the
arc length $V(\theta(t)) = \max_{i,j \in \until n} | \theta_{i}(t) -
\theta_{j}(t)|$ is a Lyapunov function candidate for phase
synchronization. Indeed, it can be shown that $V(\theta(t))$ decreases
along trajectories of \eqref{eq: coupled oscillator model - phase sync - 2}
for $\theta(0) \in \bararc(\gamma)$ and for all $\gamma < \pi$.  The
analysis is complicated by the following fact: the function $V(\theta(t))$
is continuous but not necessarily differentiable when the maximum geodesic
distance (that is, the right-hand side of \eqref{eq: contraction Lyapunov
  function}), is attained by more than one pair of oscillators. We omit the
explicit calculations here and refer to
\cite{ZL-BF-MM:07,UM-AP-FA:09,SYH-TH-KJH:10,FD-FB:09z,FD-FB:10w} for a detailed analysis.

{\em Statement 3):} By statement 2), the set $\arc(\pi)$ is positively
invariant, and for $\theta(0) \in \arc(\pi)$ the average $\sum_{i=1}^{n}
\theta_{i}(t)/n$ is well defined for $t \geq 0$.  A summation over all
equations of the model \eqref{eq: coupled oscillator model - phase sync}
yields $\sum\nolimits_{i=1}^{n} \dot \theta_{i}(t) = 0$, or equivalently,
$\sum_{i=1}^{n} \theta_{i}(t)$ is constant for all $t \geq 0$.  In
particular, for $t=0$ we have that $\sum_{i=1}^{n} \theta_{i}(t) =
\sum_{i=1}^{n} \theta_{i}(0)$ and for a phase-synchronized solution we have
that $\sum_{i=1}^{n} \subscr{\theta}{sync} = \sum_{i=1}^{n}
\theta_{i}(0)$. Hence, the explicit synchronization phase is given by
$\sum_{i=1}^{n} \theta_{i}(0)/n$. In the original coordinates (non-rotating
frame) the synchronization phase is given by $\sum_{i=1}^{n}
\theta_{i}(0)/n + \omega t \pmod{2\pi}$.

{\em Statement 4):} 
Given the invariance of the set $\bararc(\gamma)$ for any $\gamma<\pi$,
the system \eqref{eq: coupled oscillator model - phase sync - 2} can be
analyzed as a linear time-varying consensus system with initial condition
$\theta(0) \in \bararc(\gamma)$, and bounded time-varying weights
$b_{ij}(\theta(t)) \in {[ a_{ij} \sinc(\gamma),a_{ij} ]}$ for all $t \geq
0$. The worst-case convergence rate $\subscr{\lambda}{ps}$ can then be
obtained by a standard symmetric consensus analysis, see
\cite{NC-MWS:08,AJ-NM-MB:04,FD-FB:09z,FD-FB:10w}. For instance, it can be
shown that the deviation of the angles $\theta(t)$ from their average,
$\|\theta(t) - (\sum_{i=1}^{n} \theta(t)/n) \fvec 1_{n} \|_{2}^{2}$ (the
{\em disagreement function}) decays exponentially with rate
$\subscr{\lambda}{ps}$.

{\em Statement 5):} By statement 1), all solutions of system \eqref{eq:
  coupled oscillator model - phase sync} converge to the set of equilibria,
which equals the set of critical points of the potential $U(\theta)$. By
the definition of $\mathbb S^{1}$-synchronizing graphs, the
phase-synchronized equilibrium manifold $\bararc(0)$ is the only stable
equilibrium set, and all others are unstable. Hence, for all initial
condition $\theta(0) \in \torus^{n}$, which are not on the stable
manifolds of unstable equilibria, the corresponding solution $\theta(t)$
will reach the phase-synchronized equilibrium manifold $\bararc(0)$.
\end{proof}

\begin{remark}{(\bf Control-theoretic perspective on synchronization)}
As established in Theorem \ref{Theorem: Phase sync}, the set of
phase-synchronized solutions $\bararc(0)$ of the coupled oscillator
model \eqref{eq: coupled oscillator model} is locally stable provided that
all natural frequencies are identical. For non-uniform (but sufficiently
identical) natural frequencies, phase synchronization is not possible but a
certain degree of phase cohesiveness can still be achieved.
Hence, the coupled oscillator model \eqref{eq: coupled oscillator model}
can be regarded as an exponentially stable system subject to the
disturbance $\omega \in \fvec 1_{n}^{\perp}$, and classic control-theoretic
concepts such as input-to-state stability, practical stability, and
ultimate boundedness \cite{HKK:02} or their incremental versions
\cite{DA:02} can be used to study synchronization. In control-theoretic
terminology, synchronization and phase cohesiveness can then also be
described as ``practical phase synchronization''.
Compared to prototypical nonlinear control examples, various additional
challenges arise in the analysis of the coupled oscillator model \eqref{eq:
  coupled oscillator model} due to the bounded and non-monotone sinusoidal
coupling and\,the compact state space $\torus^{n}$ containing numerous
equilibria; see the analysis approaches in Section \ref{Section: Analysis
  of Synchronization} and \cite{AF-AC-WPL:11,FD-FB:10w,FD-FB:09z}.
\oprocend
\end{remark}

\subsection{Phase Balancing}

In general, only few results are known about the phase balancing
problem. This asymmetry is partially caused by the fact that phase
synchrony is required in more applications than phase balancing. Moreover,
the phase-synchronized set $\bararc(0)$ admits a very simple geometric
characterization, whereas the phase-balanced set $\bal$ has a complicated
structure consisting of numerous disjoint subsets. The number of these
subsets grows with the number of nodes $n$ in a combinatorial fashion.

Consider the coupled oscillator model \eqref{eq: coupled oscillator model - phase sync} with identical natural frequencies. By inverting the direction of time,\,we~get
\begin{equation}
	\dot \theta_{i}
	=
	\sum\nolimits_{j=1}^{n} a_{ij} \sin(\theta_{i}-\theta_{j})
	\,, \qquad i \in \until n
	\label{eq: coupled oscillator model - phase balance}
	\,.
\end{equation}
In the following, we say that the interaction graph $G(\mc V,\mc E,A)$ is {\em circulant} if the adjacency matrix $A=A^{T}$ is a circulant matrix. Circulant graphs are highly symmetric graphs including complete graphs, bipartite graphs, and ring graphs.%
\footnote{Further info on circulant graphs and a gallery can be found at {\tt http://mathworld.wolfram.com/CirculantGraph.html}.}
For circulant and uniformly weighted graphs, the coupled oscillator model \eqref{eq: coupled oscillator model - phase balance} achieves phase balancing.
The following theorem summarizes different results, which were originally presented in \cite{RS-DP-NEL:07,RS-DP-NEL:08}.

\begin{theorem}{\bf(Phase balancing)}
\label{Theorem: Phase balancing}%
Consider the coupled oscillator model \eqref{eq: coupled oscillator model - phase balance}
with a connected, uniformly weighted, and circulant graph $G(\mc V,\mc E,A)$.
The following statements hold:
\begin{enumerate}

	\item[1)] {\bf Local phase balancing:} The phase-balanced set $\bal$ is locally asymptotically stable; and

	\item[2)] {\bf Almost global stability:} If the graph $G(\mc V,\mc E,A)$ is complete, then
	the region of attraction of the stable
  phase-balanced set $\bal$ is almost all of
  $\torus^{n}$.


\end{enumerate}
\end{theorem}

The proof of Theorem \ref{Theorem: Phase balancing} follows a similar reasoning as the proof of Theorem \ref{Theorem: Phase sync}: convergence is established by potential function arguments and local (in)stability of equilibria by Jacobian arguments. We omit the proof here and refer to \cite[Theorem 1]{RS-DP-NEL:07} and \cite[Theorem 2]{RS-DP-NEL:08} for details.

For general connected graphs, the conclusions of Theorem \ref{Theorem: Phase balancing} are not true. As a remedy to achieve locally stable and globally attractive phase balancing, higher order models need to be considered, see the models proposed in \cite{LS-AS-RS:06,RS-DP-NEL:08}.

\subsection{Synchronization in Complete Networks}

For a complete coupling graph with uniform weights $a_{ij} = K/n$, where
$K>0$ is the {\em coupling gain}, the coupled oscillator model \eqref{eq:
  coupled oscillator model} reduces to the celebrated Kuramoto~model
\begin{equation}
	\dot \theta_{i}
	=
	\omega_{i} - \frac{K}{n} \sum\nolimits_{j=1}^{n} \sin(\theta_{i}-\theta_{j})
	\,, \qquad i \in \until n
	\label{eq: Kuramoto model}
	\,.
\end{equation}
By means of the order parameter $r e^{\mathrm{i} \psi} = \frac{1}{n} \sum_{j=1}^{n} e^{\mathrm{i} \theta_{j}}$, the Kuramoto model \eqref{eq:
  Kuramoto model} can be rewritten in the insightful\,form
\begin{equation}
	\dot \theta_{i} = \omega_{i} - K r \sin(\theta_{i} - \psi)
	\,,
	\quad i \in \until n
	\label{eq: Kuramoto model - order parameter}
	\,.
\end{equation}
Equation \eqref{eq: Kuramoto model - order parameter} gives the intuition
that the oscillators synchronize by coupling to a mean field represented by
the order parameter $r e^{\mathrm{i} \psi}$.  Intuitively, for small
coupling strength $K$ each oscillator rotates with its natural frequency
$\omega_{i}$, whereas for large coupling strength $K$ all angles
$\theta_{i}(t)$ will be entrained by the mean field $r e^{\mathrm{i} \psi}$
and the oscillators synchronize.
The threshold from incoherence to synchronization occurs for some critical
coupling $\subscr{K}{critical}$. This phase transition has been the source
of numerous investigations starting with Kuramoto's analysis
\cite{YK:75,YK:84}.
Various necessary, sufficient, implicit, and explicit estimates of the critical coupling strength $\subscr{K}{critical}$ for both the on-set as well as the ultimate stage of synchronization have been proposed \cite{YK:75,FD-FB:10w,SHS:00,JAA-LLB-CJPV-FR-RS:05,RS-DP-NEL:07,FD-FB:09z,YK:84,GBE:85,EAM-EB-SHS-PS-TMA:09,JLvH-WFW:93,NC-MWS:08,AJ-NM-MB:04,SJC-JJS:10b,GSS-UM-FA:09,FDS-DA:07,MV-OM:08,REM-SHS:05,DA-JAR:04,EC-PM:09,PM:06,MV-OM:09,AS:09,LB-LS-ADG:09,AF-AC-WPL:11,SYH-TH-KJH:10,SYH-MS:11}, and we refer to \cite{FD-FB:10w} for a comprehensive\,overview.

The mean field approach to the equations \eqref{eq: Kuramoto model - order parameter} can be made mathematically rigorous by a time-scale separation \cite{SYH-MS:11} or in the continuum limit as the number of oscillators tends to infinity and the natural frequencies $\omega$ are sampled from a distribution function $g: \real \to \real_{\geq 0}$. 
%
In the continuum limit and for a symmetric, continuous, and unimodal  distribution $g(\omega)$, Kuramoto himself showed in an insightful and ingenuous analysis \cite{YK:75,YK:84} that the incoherent state (a uniform distribution of the oscillators on the unit circle $\mycircle$) supercritically bifurcates for the critical coupling strength 
\begin{equation}
	\subscr{K}{critical} = \frac{2}{\pi g(0)}
	\label{eq: Kuramoto bound - continuum limit}
	\,.
\end{equation}
In \cite{GBE:85,JLvH-WFW:93,JAA-LLB-CJPV-FR-RS:05}, it was found that the bipolar (bimodal double-delta) distribution (respectively the uniform distribution) yield the largest (respectively smallest) threshold $\subscr{K}{critical}$ over all distributions $g(\omega)$ with bounded support. We refer \cite{SHS:00,JAA-LLB-CJPV-FR-RS:05} for further references and to \cite{YH-PGM-SPM-UVS:12,CH:10,VM-OM:11,EAM-EB-SHS-PS-TMA:09} for recent contributions on the continuum limit model. 


In the finite-dimensional case, the necessary synchronization condition \eqref{eq: necessary sync condition - 2} gives a lower bound for
$\subscr{K}{critical}$ as 
\begin{equation}
	K \geq \frac{n}{2(n-1)} \cdot ( \subscr{\omega}{max} - \subscr{\omega}{min} )
	\,.
	\label{eq: Kuramoto bound - necessary}
\end{equation}


Three recent articles \cite{MV-OM:08,REM-SHS:05,DA-JAR:04} independently derived a set of implicit consistency equations for the {\it exact} critical coupling strength $\subscr{K}{critical}$ for which synchronized solutions exist. Verwoerd and Mason provided the following implicit formulae to compute $\subscr{K}{critical}$ \cite[Theorem 3]{MV-OM:08}:
\begin{align}
	\label{eq: Kuramoto bound - exact}
\begin{split}
	&
	\subscr{K}{critical} 
	=
	n u^{*} / \sum\nolimits_{i=1}^{n} \sqrt{1 - (\Omega_{i}/u^{*})^{2}}
	\,,
	\\
	&\!
	2 \sum\nolimits_{i=1}^{n} \!\sqrt{1- (\Omega_{i}/u^{*})^{2}} = \sum\nolimits_{i=1}^{n} \! 1/\sqrt{1 - (\Omega_{i}/u^{*})^{2}}
	,
	\end{split}
\end{align}
where $\Omega_{i} = \omega_{i} - \subscr{\omega}{sync}$ and $u^{*} \in
[\norm{\Omega}_{\infty},2 \norm{\Omega}_{\infty}]$. The~implicit formulae
\eqref{eq: Kuramoto bound - exact} can also be extended to bipartite graphs
\cite{MV-OM:09}. A local stability analysis is carried out in \cite{REM-SHS:05,DA-JAR:04}.

From the point of analyzing or designing a sufficiently strong coupling,
the exact formulae \eqref{eq: Kuramoto bound - exact} have three
drawbacks. First, they are implicit and thus not suited for performance or
robustness estimates in case of additional coupling strength for a given
$K>\subscr{K}{critical}$. Second, the corresponding region of attraction of
a synchronized solution is unknown. Third and finally, the particular
natural frequencies $\omega_{i}$ are typically time-varying, uncertain, or
even unknown in the applications listed in Section \ref{Section:
  Introduction}. In this case, the exact~value of $\subscr{K}{critical}$
needs to be estimated in continuous time, or a conservatively strong
coupling $K \!\gg\! \subscr{K}{critical}$ has to be~chosen.

The following theorem states an explicit bound on the critical coupling strength together with performance estimates, convergence rates, and a guaranteed semi-global region of attraction for synchronization. This bound is tight and thus necessary and sufficient when considering arbitrary distributions of the natural frequencies with compact support. The result has been originally presented in \cite[Theorem 4.1]{FD-FB:10w}.

\begin{theorem}{\bf(Synchronization in the Kuramoto model)}
\label{Theorem: Kuramoto frequency sync}
Consider the Kuramoto model \eqref{eq: Kuramoto model} with natural
frequencies $\omega = (\omega_1,\dots,\omega_n)$ and coupling strength $K$.
The following three statements are equivalent:

\begin{enumerate}
\item[(i)] the coupling strength $K$ is larger than the maximum
  non-uniformity among the natural frequencies, that is,
  \begin{equation}
    K > \subscr{K}{critical} \triangleq 
    \subscr{\omega}{max} - \subscr{\omega}{min}
    \label{eq: Kuramoto bound - sufficient}
    \;;
  \end{equation}

\item[(ii)] there exists an arc length $\subscr{\gamma}{max}\in
  {]\pi/2,\pi]}$ such that the Kuramoto model \eqref{eq: Kuramoto model}
synchronizes exponentially for all possible distributions of the natural
frequencies $\omega_i$ supported on the compact interval
$[\subscr{\omega}{min},\subscr{\omega}{max}]$ and for all initial phases
$\theta(0) \in \arc(\subscr{\gamma}{max})$; and

\item[(iii)] there exists an arc length $\subscr{\gamma}{min} \in
  {[0,\pi/2[}$ such that the Kuramoto model \eqref{eq: Kuramoto model} has
      a locally exponentially stable synchronization manifold in
      $\bararc(\subscr{\gamma}{min})$ for all possible distributions of
      the natural frequencies $\omega_{i}$ supported on the compact
      interval $[\subscr{\omega}{min},\subscr{\omega}{max}]$.

\end{enumerate}

\smallskip 

If the three equivalent conditions (i), (ii), and (iii) hold, then the
ratio $\subscr{K}{critical}/K$ and the arc lengths $\subscr{\gamma}{min}
\in {[0,\pi/2[}$ and $\subscr{\gamma}{max} \in {]\pi/2,\pi]}$ are related
uniquely via $\sin(\subscr{\gamma}{min}) = \sin(\subscr{\gamma}{max}) =
{\subscr{K}{critical}}/K$, and the following statements hold: \smallskip

\begin{enumerate}

\item[1)] {\bf phase cohesiveness:} the set $\bararc(\gamma) \subseteq \bar\Delta_{G}(\gamma)$ is positively invariant for every $\gamma \in [\subscr{\gamma}{min},\subscr{\gamma}{max}]$,\,and each trajectory starting in $\arc(\subscr{\gamma}{max})$ approaches asymptotically $\bararc(\subscr{\gamma}{min})$;
	
	  \item[2)] {\bf frequency synchronization:} the asymptotic synchronization frequency is
  the average frequency $\subscr{\omega}{sync} = \frac{1}{n}
  \sum_{i=1}^{n} \omega_{i}$, and, given phase cohesiveness in
  $\bararc(\gamma)$ for some fixed $\gamma <\pi/2$, the exponential
  synchronization rate is no worse than $\subscr{\lambda}{K} = - K
  \cos(\gamma)$; and
  
\item[3)] {\bf order parameter:} the asymptotic value of the magnitude of
  the order parameter, denoted by $r_{\infty} \triangleq \lim_{t \to
    \infty} \frac{1}{n} |\sum_{j=1}^{n} e^{\mathrm{i}\theta_{j}(t)} |$, is
  bounded as
  \begin{equation*}
    1  \geq r_{\infty} \geq \cos\!\left(\frac{\subscr{\gamma}{min}}{2}\right) = \sqrt{\frac{1+\sqrt{1- (\subscr{K}{critical}/K)^{2}}}{2}}
    \,.
  \end{equation*}  

\end{enumerate}
\end{theorem}

\smallskip\begin{proof}
In the following, we sketch the proof of Theorem \ref{Theorem: Kuramoto frequency sync} and refer to \cite[Theorem 4.1]{FD-FB:10w} for further details. 

{\em Implication (i) $\implies$ (ii):} In a first step, it is shown that
the phase cohesive set $\bararc(\gamma)$ is positively invariant for
every $\gamma \in [\subscr{\gamma}{min},\subscr{\gamma}{max}]$. By
assumption, the angles $\theta_i(t)$ belong to the set $\bararc(\gamma)$
at time $t=0$, that is, they are all contained in an arc of length
$\gamma$. We aim to show that all angles remain in $\bararc(\gamma)$ for
all subsequent times $t>0$ by means of the contraction Lyapunov function
\eqref{eq: contraction Lyapunov function}.
Note that $\bararc(\gamma)$ is positively invariant if and only if
$V(\theta(t))$ does not increase at any time $t$ such that $V(\theta(t))=
\gamma$. The {\it upper Dini derivative} of $V(\theta(t))$ along
trajectories of \eqref{eq: Kuramoto model} is given by
\begin{equation*}
	D^{+} V (\theta(t))
	=
	\lim_{h \downarrow 0}\sup \frac{V(\theta(t+h)) - V(\theta(t))}{h}
	\,.
\end{equation*}
Written out in components and after trigonometric simplifications
\cite{FD-FB:10w}, we obtain that the derivative is bounded as
\begin{equation*}
	D^{+} V (\theta(t))
	\leq
	 \subscr{\omega}{max} - \subscr{\omega}{min} - K \sin(\gamma)
	\,.
\end{equation*}
It follows that the length of the arc formed by the angles is
non-increasing in $\bararc(\gamma)$ if and only if
\begin{equation}
	K \sin(\gamma)
	\geq
	\subscr{K}{critical}
	\label{eq: critical bound}
	\,,
\end{equation}
where $\subscr{K}{critical}$ is as stated in equation \eqref{eq: Kuramoto
  bound - sufficient}. For $\gamma \in {[0,\pi]}$ the left-hand side of
\eqref{eq: critical bound} is a concave function of $\gamma$ that achieves
its maximum at $\gamma^{*}=\pi/2$. Therefore, there exists an open set of
arc lengths $\gamma \in {[0,\pi]}$ satisfying equation \eqref{eq: critical
  bound} if and only if equation \eqref{eq: critical bound} is true with
the strict equality sign at $\gamma^{*}=\pi/2$, which corresponds to
condition \eqref{eq: Kuramoto bound - sufficient}. Additionally, if these
two equivalent statements are true, then there exists a unique
$\subscr{\gamma}{min}\in {[0,\pi/2[}$ and a $\subscr{\gamma}{max}\in
    {]\pi/2,\pi]}$ that satisfy equation \eqref{eq: critical bound} with
the equality sign, namely $\sin(\subscr{\gamma}{min}) =
\sin(\subscr{\gamma}{max}) = {\subscr{K}{critical}}/K$. For every $\gamma
\in {[\subscr{\gamma}{min},\subscr{\gamma}{max}]}$ it follows that the
arc-length $V(\theta(t))$ is non-increasing, and it is strictly decreasing
for $\gamma \in {]\subscr{\gamma}{min},\subscr{\gamma}{max}[}$. Among other
things, this shows that statement (i) implies statement 1). By means of
Lemma \ref{Lemma: phase cohesiveness and order parameter}, statement 3)
then follows from statement 1).

The frequency dynamics of the Kuramoto model \eqref{eq: Kuramoto model} can be obtained by differentiating the Kuramoto model \eqref{eq: Kuramoto model} as
\begin{equation}
	\dt \dot{\theta_{i}}
	=
	\sum\nolimits_{j=1}^{n} \tilde a_{ij}(t) \,(\dot{\theta}_{j} - \dot{\theta}_{i})
	\label{eq: consensus system for dot theta}
	\,,
\end{equation}
where $\tilde a_{ij}(t) = (K/n) \cos(\theta_{i}(t)-\theta_{j}(t))$. For $K > \subscr{K}{critical}$, we just proved that for every $\theta(0) \in \arc(\subscr{\gamma}{max})$ and for all $\gamma \in {]\subscr{\gamma}{min},\subscr{\gamma}{max}[}$ there exists a finite time $T \geq 0$ such that $\theta(t) \in \bararc(\gamma)$ for all $t \geq T$. Consequently, the terms $\tilde a_{ij}(t)$ are strictly positive for all $t \geq T$. Notice also that system \eqref{eq: consensus system for dot theta} evolves on the tangent space of $\torus^{n}$, that is, the Euclidean space $\real^{n}$. Now fix $\gamma \in {]\subscr{\gamma}{min},\pi/2[}$ and let $T \geq 0$ such that $\tilde a_{ij}(t)> 0$ for all $t \geq T$. In this case, the frequency dynamics \eqref{eq: consensus system for dot theta} can be analyzed as linear time-varying consensus system.
Consider the {\it disagreement vector} $ x = \dot{\theta} -
\subscr{\omega}{sync} \fvec 1_{n} $ as an error coordinate. By standard
consensus arguments \cite{ROS-JAF-RMM:07,WR-RWB-EMA:07,FB-JC-SM:09}, it can
be shown that the disagreement vector satisfies $ \|x(t)\| \!\leq\!
\|x(0)\| e^{- \subscr{\lambda}{K} t} $ for all\,$t \geq T$. This proves
statement 2) and the implication (i)\,$\implies$\,(ii).

{\em Implication (ii) $\implies$ (i):}
To show that condition \eqref{eq: Kuramoto bound - sufficient} is also necessary for synchronization, it suffices to construct a counter example for which $K \leq \subscr{K}{critical}$ and the oscillators do not achieve exponential synchronization even though all $\omega_{i} \in {[\subscr{\omega}{min},\subscr{\omega}{max}]}$ and $\theta(0) \in \arc(\gamma)$ for every $\gamma \in {]\pi/2,\pi]}$. 
 A basic instability mechanism under which synchronization breaks down is caused by a bipolar distribution of the natural frequencies. Let the index set $\until n$ be partitioned by the two non-empty sets $\mc I_{1}$ and $\mc I_{2}$. Let $\omega_{i} = \subscr{\omega}{min}$ for $i \in \mc I_{1}$ and $\omega_{i} = \subscr{\omega}{max}$ for $i \in \mc I_{2}$, and assume that at some time $t \geq 0$ it holds that $\theta_{i}(t) \!=\! - \gamma/2$ for $i \in \mc I_{1}$ and $\theta_{i}(t) \!=\! + \gamma/2$ for $i \in \mc I_{2}$ and for some $\gamma \in {[0,\pi[}$. By construction, at time $t$ all oscillators are contained in an arc of length $\gamma \in {[0,\pi[}$. Assume now that $K\!<\!\subscr{K}{critical}$ and the oscillators synchronize. It can be shown \cite{FD-FB:10w} that the evolution of the arc length $V(\theta(t))$ satisfies the {\em equality}
\begin{equation}
	D^{+} V (\theta(t))
	=
	\subscr{\omega}{max} - \subscr{\omega}{min} - K \sin(\gamma)
	\label{eq: critical bound - necessary}
	\,.
\end{equation}
Clearly, for $K < \subscr{K}{critical}$ the arc length $V(\theta(t)) = \gamma$ is increasing for any arbitrary $\gamma \in {[0,\pi]}$. Thus, the phases are not bounded in $\bararc(\gamma)$. This contradicts the assumption that the oscillators synchronize for $K < \subscr{K}{critical}$ from every initial condition $\theta(0) \in \bararc(\gamma)$. For $K = \subscr{K}{critical}$, we know from \cite{REM-SHS:05,DA-JAR:04} that phase-locked equilibria have a zero eigenvalue with a two-dimensional Jacobian block, and thus synchronization cannot occur. This instability via a two-dimensional Jordan block is also visible in \eqref{eq: critical bound - necessary} since $D^{+} V (\theta(t))$ is increasing for $\theta(t) \in \bararc(\gamma)$, $\gamma \in {]\pi/2,\pi]}$ until all oscillators change orientation, just as in the example in Subsection \ref{Subsection: 2d example}. This proves the implication (ii)\,$\implies$\,(i).

{\em Equivalence (i),(ii) $\Leftrightarrow$ (iii):} The proof relies on Jacobian arguments and will be omitted here, see \cite{FD-FB:10w} for details.
\end{proof}

Theorem \ref{Theorem: Kuramoto frequency sync} places a hard bound on the critical coupling strength $\subscr{K}{critical}$ for all distributions of $\omega_{i}$ supported on the compact interval $[\subscr{\omega}{min},\subscr{\omega}{max}]$. For a particular distribution $g(\omega)$ supported on $[\subscr{\omega}{min},\subscr{\omega}{max}]$ the bound \eqref{eq: Kuramoto bound - sufficient} is only sufficient and possibly a factor 2 larger than the necessary bound \eqref{eq: Kuramoto bound - necessary}. The exact critical coupling lies somewhere in between and can be obtained from the implicit equations \eqref{eq: Kuramoto bound - exact}. 

\begin{figure}[hbt]
	\centering{
	\includegraphics[width=0.99\columnwidth]{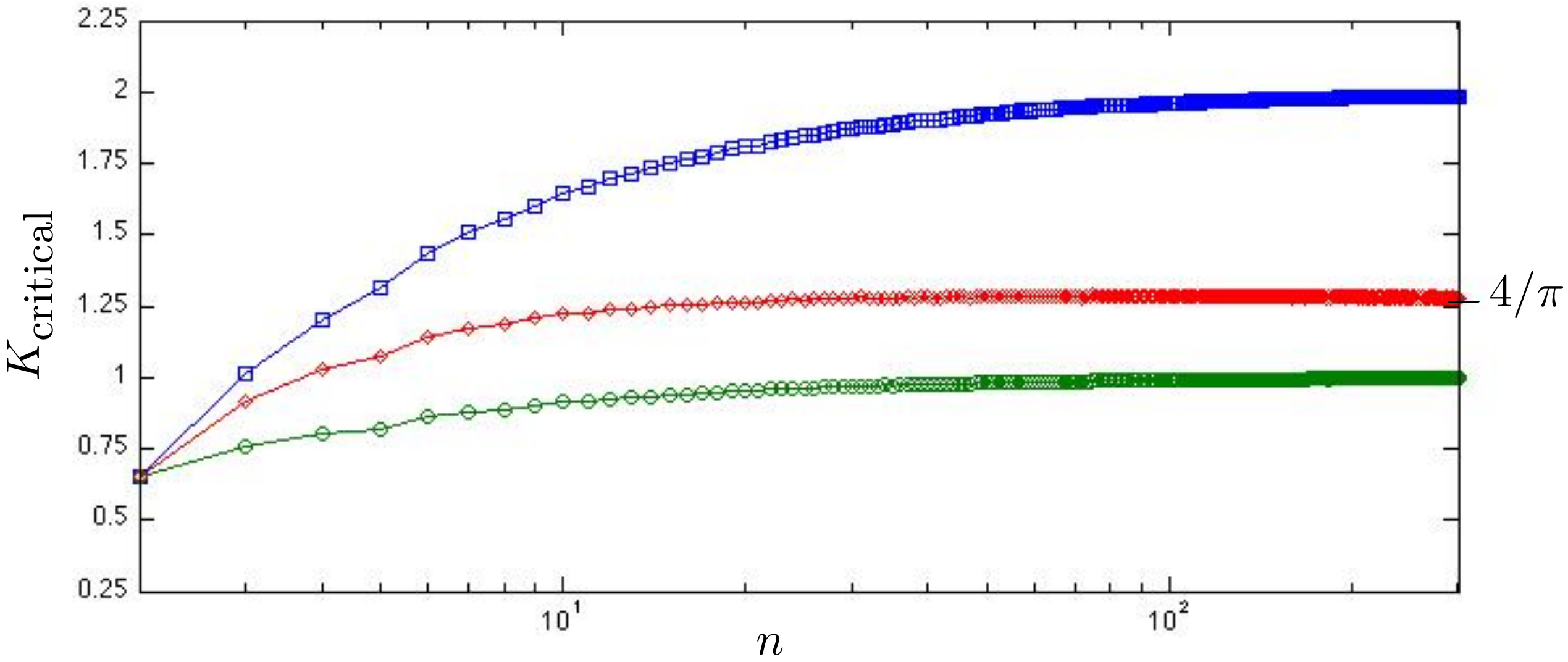}
	\caption{Statistical analysis of the necessary and explicit bound \eqref{eq: Kuramoto bound - necessary} ({\bf\footnotesize$\color{green}\lozenge$}), the exact and implicit bound \eqref{eq: Kuramoto bound - exact} ({\bf\large$\color{red}{\circ}$}), and the sufficient, tight, and explicit bound \eqref{eq: Kuramoto bound - sufficient} ({\bf\tiny$\color{blue}{\square}$}) for $n \in [2,300]$ oscillators in a semi-log plot, where the coupling gains for each $n$ are averaged over 1000 simulations.}
	\label{Fig: Kuramoto bounds for uniform distribution}
	}
\end{figure}

Since the bound \eqref{eq: Kuramoto bound - sufficient} on $\subscr{K}{critical}$ is exact \cite{FD-FB:10w} for the worst-case bipolar distribution $\omega_{i} \in \{\subscr{\omega}{min},\subscr{\omega}{max}\}$, Figure \ref{Fig: Kuramoto bounds for uniform distribution} reports numerical findings for the other extreme case \cite{GBE:85} of a uniform distribution $g(\omega) = 1/2$ supported for $\omega_{i} \in [-1,1]$. 
All three displayed bounds are identical for $n=2$ oscillators. As $n$
increases, the sufficient bound \eqref{eq: Kuramoto bound - sufficient}
converges to the\,width $\subscr{\omega}{max}-\subscr{\omega}{min}=2$ of
the support of $g(\omega)$, and the necessary bound \eqref{eq: Kuramoto
  bound - necessary} accordingly to half of the width. The exact bound
\eqref{eq: Kuramoto bound - exact} converges to
$4(\subscr{\omega}{max}-\subscr{\omega}{min})/(2\pi) \!=\! 4/\pi$ in
agreement with condition \eqref{eq: Kuramoto bound - continuum limit}
predicted for the continuum limit.

Finally, let us mention that Theorem \ref{Theorem: Kuramoto frequency sync}
can be extended to discontinuously switching and slowly time-varying
natural frequencies \cite{FD-FB:10w}. For a particular sampling distribution 
$g(\omega)$, the critical quantity in condition \eqref{eq: Kuramoto bound - sufficient}, the support $\subscr{\omega}{max} - \subscr{\omega}{min}$, can be estimated by extreme value statistics, see \cite{JCB-LDV-MJP:11}.

\subsection{Synchronization in Sparse Networks}

As summarized in Subsection \ref{Subsection: Synchronization Conditions},
the quest for sharp and concise synchronization for non-complete coupling
graph $G(\mc V,\mc E,A)$ is an important and outstanding problem emphasized
in every review article on coupled oscillator networks
\cite{SHS:01,JAA-LLB-CJPV-FR-RS:05,FD-FB:10w,SHS:00,AA-ADG-JK-YM-CZ:08,SB-VL-YM-MC-DUH:06}. The
approaches known for phase synchronization in arbitrary graphs or the
contraction approach to frequency synchronization (used in the proof of
Theorem \ref{Theorem: Kuramoto frequency sync}) do not generally extend to
arbitrary natural frequencies $\omega \in \fvec 1_{n}^{\perp}$ and connected
coupling graphs $G(\mc V,\mc E,A)$, or do so only under extremely
conservative conditions.

One Lyapunov function advocated for classic Kuramoto oscillators \eqref{eq:
  Kuramoto model} is the function $W:\, \arc(\pi) \to \real$ defined for
angles $\theta$ in an open semi-circle and given by
\cite{NC-MWS:08,AJ-NM-MB:04}
\begin{equation}
	W(\theta) = \frac{1}{4} 
        \sum\nolimits_{i,j=1}^{n} | \theta_{i} - \theta_{j}|^{2}
	= \frac{1}{2} \left\| B_{c}^{T} \theta \right\|_{2}^{2}
	\label{eq: two-norm Lyapunov function}
	\,,
\end{equation}
where $B_{c} \in \real^{n \times (n(n-1)/2)}$ is an incidence matrix of the complete graph. As shown in \cite[Theorem 4.4]{FD-FB:09z}, the Lyapunov function \eqref{eq: two-norm Lyapunov function} generalizes also to the coupled oscillator model \eqref{eq: coupled oscillator model}. Indeed, an even more general model\,is considered in \cite{FD-FB:09z}, and a Lyapunov analysis yields the following result. 

\begin{theorem}
\label{Theorem: Synchronization in coupled oscillators I}
{\bf(Frequency synchronization I)}
Consider the coupled oscillator model \eqref{eq: coupled oscillator model} with a connected graph $G(\mc V,\mc E,A)$ and $\omega \in \fvec 1_{n}^{\perp}$. Assume that the algebraic connectivity is larger than a critical value, that is,
\begin{equation}
  \label{eq: frequ sync condition - 1}
  \lambda_{2}(L) 
  > \subscr{\lambda}{critical} 
  \triangleq 
  \norm{B_{c}^{T} \omega}_{2}
  \,,
\end{equation}
where $B_{c} \in \real^{n \times n(n-1)/2}$ is the incidence matrix of the complete graph.
Accordingly, define $\subscr{\gamma}{max} \in {]\pi/2,\pi]}$ and
 $\subscr{\gamma}{min} \in {[0,\pi/2[}$ as unique solutions to $(\pi/2)
     \cdot \sinc(\subscr{\gamma}{max}) \!=\!  \sin(\subscr{\gamma}{min})
     \!=\!  \subscr{\lambda}{critical} / \lambda_{2}(L)$. The following
     statements hold:
\begin{enumerate}	

\item[1)] {\bf phase cohesiveness:} the set $\bigl\{\theta \in
  \arc(\pi):\, \|B_{c}^{T}\theta\|_{2} \leq \gamma \bigr\} \subseteq \bar\Delta_{G}(\gamma)$ is positively
  invariant for every $\gamma \in
  [\subscr{\gamma}{min},\subscr{\gamma}{max}]$, and each trajectory
  starting in the set $\bigl\{\theta \in \arc(\pi):\, \norm{B_{c}^{T}\theta}_{2}
  < \subscr{\gamma}{max}\bigr\}$ asymptotically reaches the set $\bigl\{\theta \in
  \arc(\pi):\,\|B_{c}^{T}\theta\|_{2} \leq \subscr{\gamma}{min} \bigr\}$;
  and
  
\item[2)] {\bf frequency synchronization:} for every $\theta(0) \in
  \arc(\pi)$ with $\norm{B_{c}^{T}\theta(0)}_{2} < \subscr{\gamma}{max}$
  the frequencies $\dot{\theta}_{i}(t)$ synchronize exponentially to the
  average frequency $\subscr{\omega}{sync} = \frac{1}{n} \sum_{i=1}^{n}
  \omega_{i}$, and, given phase cohesiveness in $\bar\Delta_{G}(\gamma)$
  for some fixed $\gamma <\pi/2$, the exponential synchronization rate is
  no worse than $\subscr{\lambda}{fe} = - \lambda_{2}(L) \cos(\gamma)$.
  
\end{enumerate}
\end{theorem}
\smallskip

The proof of Theorem \ref{Theorem: Synchronization in coupled oscillators I} follows a similar ultimate-boundedness strategy as the proof of Theorem \ref{Theorem: Kuramoto frequency sync} by using the Lyapunov function \eqref{eq: two-norm Lyapunov function}. It can be found in Appendix \ref{Subsection: Proof of freq sync in sparse graphs I}.

For classic Kuramoto oscillators \eqref{eq: Kuramoto model}, condition \eqref{eq: frequ sync condition - 1} reduces to 
  $K > \norm{B_c^{T}\omega}_{2}$. Clearly, the condition $K \!>\! \norm{B_c^{T}\omega}_{2}$ is more conservative than the bound \eqref{eq: Kuramoto bound - sufficient} which reads as $K > \|B_{c}^{T} \omega \|_{\infty} = \subscr{\omega}{max} - \subscr{\omega}{min}$.
One reason for this conservatism is that the analysis leading to condition \eqref{eq: frequ sync condition - 1} requires {\em all} phase distances $|\theta_{i} - \theta_{j}|$ to be bounded, whereas according to Lemma \ref{Lemma: stable equilibrium in pi/2 arc} only {\em pairwise}\,\,phase distances $|\theta_{i} - \theta_{j}|$, $\{i,j\} \in \mc E$, need to be bounded for stable synchronization. 
The following result exploits these weaker assumptions and states a sharper (but only local) synchronization condition. 
%

\begin{theorem}
\label{Theorem: Synchronization in coupled oscillators II}
{\bf(Frequency synchronization II)}
Consider the coupled oscillator model \eqref{eq: coupled oscillator model} with a connected graph
$G(\mc V,\mc E,A)$ and $\omega \in \fvec 1_{n}^{\perp}$. There exists a locally exponentially stable equilibrium manifold $[\theta] \in \Delta_{G}(\pi/2)$ if 
\begin{equation}
	\lambda_{2}(L) > \left\| B^{T} \omega \right\|_{2}
	\label{eq: frequ sync condition - 2}
	\,.
\end{equation}
%
Moreover, if condition \eqref{eq: frequ sync condition - 2} holds, then $[\theta]$ is phase cohesive in $\{\theta \in \torus^{n}:\, \| B^{T} \theta \|_{2} \leq \subscr{\gamma}{min} \} \subseteq \bar\Delta_{G}(\subscr{\gamma}{min})$, where $\subscr{\gamma}{min} \in {[0,\pi/2[}$ satisfies 
$\sin(\subscr{\gamma}{min}) 
=
\|  B^{T}\omega \|_{2} / \lambda_{2}(L)$. 
\end{theorem}
\smallskip

The strategy to prove Theorem \ref{Theorem: Synchronization in coupled oscillators II} is inspired by the ingenuous analysis in \cite[Section IIV.B]{AJ-NM-MB:04}. 
It relies on the insight gained from Lemma \ref{Lemma: stable equilibrium in pi/2 arc} that any synchronization manifold $[\theta] \in \Delta_{G}(\pi/2)$ is locally stable, and it formulates the existence of such a synchronization manifold as a fixed point problem.
Here, we follow the basic proof strategy in \cite{AJ-NM-MB:04}, but we provide a more accurate result together with a self-contained proof which is reported in Appendix \ref{Subsection: Proof of freq sync in sparse graphs II}.



\section{Conclusions and Open Research Directions}
\label{Section: Conclusions}

In this paper we introduced the reader to the coupled oscillator
model \eqref{eq: coupled oscillator model}, we reviewed several
applications, we discussed different synchronization notions, and we
presented different analysis approaches to phase synchronization,
phase balancing, and frequency
synchronization. 

Despite the vast literature, the countless applications, and the numerous
theoretic results on the synchronization properties of model~\eqref{eq:
  coupled oscillator model}, many interesting and important problems are
still open. In the following, we summarize limitations of the existing analysis approaches and present a few worthwhile directions for
future research.

First, in many applications the coupling between the oscillators is not
purely sinusoidal. For instance, phase delays in neuroscience
\cite{SMC-GBE-MCV-JMB:97}, time delays in sensor networks
\cite{OS-US-YBN-SS:08}, or transfer conductances in power networks
\cite{HDC-CCC-GC:95} lead to a ``shifted coupling'' of the form
$\sin(\theta_{i} - \theta_{j} - \varphi_{ij})$, where $\varphi_{ij} \in {[-
    \pi/2,\pi/2]}$. In this case and also for other ``skewed'' or ``symmetry-breaking'' coupling
functions, many of the presented analysis schemes either fail or lead to overly
conservative results. 
Another interesting class of oscillator networks are systems of pulse-coupled oscillators featuring hybrid dynamics: impulsive coupling at discrete time instants and uncoupled continuous dynamics otherwise. This class of oscillator networks displays a very interesting phenomenology. For instance, the behavior of identical oscillators coupled in a complete graph strongly depends on the curvature of the uncoupled dynamics \cite{AM:11}. Most of the results known for continuously-coupled oscillators still need to be extended to pulse-coupled oscillators with hybrid dynamics.

Second, in many applications
\cite{FCH-EMI:00,EB-PH-JM:03,KW-PC-SHS:98,HDC-CCC-GC:95,MB-MFS-HR-KW:02}
the coupled oscillator dynamics are not given by a simple first-order phase
model of the form \eqref{eq: coupled oscillator model}. Rather, the
dynamics are of higher order, or sometimes there is no readily available
phase variable to describe the limit cycle attracting the coupled
dynamics. The analysis of oscillator networks with more general oscillator
dynamics is largely unexplored. Whereas advances have been made for the
simple case of phase synchronization of linear or passive oscillator
networks, the case of frequency synchronization of non-identical
oscillators with higher-order dynamics is not well-studied.

Third, despite the vast scientific interest the quest for sharp, concise,
and closed-form synchronization conditions for arbitrary complex graphs has
been so far in vain
\cite{SHS:01,JAA-LLB-CJPV-FR-RS:05,FD-FB:10w,SHS:00,AA-ADG-JK-YM-CZ:08,SB-VL-YM-MC-DUH:06}. As
suggested by Lemma \ref{Lemma: Necessary sync condition}, Lemma \ref{Lemma:
  stable equilibrium in pi/2 arc}, Theorem \ref{Theorem: Kuramoto frequency
  sync}, and the proof of Theorem \ref{Theorem: Synchronization in coupled
  oscillators II}, the proper metric for the synchronization problem is the
incremental $\infty$-norm $\|B^{T} \theta\|_{\infty} = \max_{\{i,j\} \in
  \mc E} |\theta_{i} - \theta_{j}|$. In the authors' opinion, a Banach
space analysis of the coupled oscillator model \eqref{eq: coupled
  oscillator model} with the incremental $\infty$-norm will most likely
deliver the sharpest possible conditions. However, such an analysis is very
challenging for arbitrary natural frequencies $\omega \in \fvec
1_{n}^{\perp}$ and connected and weighted coupling graphs\,$G(\mc V,\mc
E,A)$. 
Recent work~\cite{FD-MC-FB:11v-arxiv} by the authors puts forth a
novel algebraic condition for synchronization with a rigorous analysis for
specific classes of graphs and with (only) a statistical validation for
generic weighted graphs.

Fourth and finally, a few interesting and open theoretical challenges
include the following.
First, most of the presented analysis approaches and conditions do not
extend to time-varying or directed coupling graphs $G(\mc V,\mc E,A)$, and alternative
methods need to be developed.
Second, most known estimates on the region of attraction of a synchronized
solution are conservative. The semi-circle estimates given in Theorem \ref{Theorem: Phase sync} and Theorem \ref{Theorem:
  Kuramoto frequency sync} rely on convexity of $\arc(\pi)$ and are overly conservative. We refer to \cite{DAW-SHS-MG:06,HDC-CCC-GC:95}
for a set of interesting results and conjectures on the region of attraction.
Third, the presented analysis approaches are restricted to synchronized
equilibria inside the set $\Delta_{G}(\pi/2)$. Other interesting
equilibrium configurations outside $\Delta_{G}(\pi/2)$ include splay state
equilibria or frequency-synchronized equilibria 
with phases spread over an entire
semi-circle.

We sincerely hope that this tutorial article stimulates further exciting
research on synchronization in coupled oscillators, both on the theoretical
side as well as in the countless applications.


\begin{appendix}

\subsection{Modeling of the spring-interconnected particles}
\label{Subsection: Modeling of the Mechanical Analog}

Consider the spring network in Figure \ref{Fig: Mechanical analog}
consisting of a group of particles constrained to rotate around a
circle of unit radius. For simplicity, we assume that the particles are allowed to move
freely on the circle and exchange their order without collisions. Each
particle is characterized by its phase angle $\theta_{i} \in \mathbb S^{1}$
and frequency $\dot \theta_{i} \in \real$, and its inertial and damping
coefficients are $M_{i}>0$ and $D_{i}>0$.

The external forces and torques acting on each particle are (i) a viscous
damping force $D_{i} \dot \theta_{i}$ opposing the direction of motion,
(ii) a non-conservative force $\omega_{i} \in \real$ along the direction of
motion depicting a preferred natural rotation frequency, and (iii) an
elastic restoring torque between interacting particles $i$ and $j$ coupled
by an ideal elastic spring with stiffness $a_{ij} > 0$ and zero rest
length.

To compute the elastic torque between the particles, we parametrize the
position of each particle $i$ by the unit vector $p_{i} = \left[
  \cos(\theta_{i}) \,,\, \sin(\theta_{i}) \right]^{T} \in \mathbb S^{1}
\subset \real^{2}$. The elastic Hookean energy stored in the springs is the
function $E:\, \torus^{n} \to \real$ given up to an additive constant by
\begin{align*}
 E(\theta) 
 &= \sum\nolimits_{\{i,j\} \in \mc E} \frac{a_{ij}}{2}  \| p_{i} - p_{j} \|_{2}^{2} \\
 &= \sum\nolimits_{\{i,j\} \in \mc E} a_{ij}  \bigl( 1 -\cos(\theta_{i})\cos(\theta_{j}) - \sin(\theta_{i})\sin(\theta_{j}) \bigr) \\
 &= \sum\nolimits_{\{i,j\} \in \mc E} a_{ij}  \bigl( 1 - \cos(\theta_{i} - \theta_{j})\bigr)
 \,,
\end{align*}
where we employed the trigonometric identity $\cos(\alpha - \beta) = \cos
\alpha \cos \beta + \sin \alpha \sin \beta$ in the last equality. Hence, we
obtain the restoring torque acting on particle $i$ as
\begin{equation*}
 T_{i}(\theta) = - \frac{\partial}{\partial\theta_{i}} \, E(\theta) = - \sum\nolimits_{\{i,j\} \in \mc E} a_{ij} \sin(\theta_{i} - \theta_{j})
 \,.
\end{equation*}
Therefore, the network of spring-interconnected particles depicted in
Figure \ref{Fig: Mechanical analog} obeys the dynamics
\begin{multline}
 M_i \ddot{\theta}_i + D_i \dot{\theta}_i = \omega_{i} - \sum\nolimits_{\{i,j\} \in \mc E} a_{ij}\sin(\theta_i-\theta_j)
 \;,\\
 i \in \until n
 \label{eq: spring network}
 \,.
\end{multline}
The coupled oscillator model \eqref{eq: coupled oscillator model} is then
obtained as the kinematic variant or the overdamped limit of the spring
network \eqref{eq: spring network} with zero inertia $M_{i} = 0$ and unit
damping $D_{i} = 1$ for all oscillators $i \in \until n$.

\subsection{Proof of Theorem \ref{Theorem: Synchronization in coupled oscillators I}}
\label{Subsection: Proof of freq sync in sparse graphs I}

Assume that $\theta(0) \in \bararc(\rho)$ for
$\rho \in {[0,\pi[}$. Recall that the angular differences are well defined
for $\theta$ in the open semi-circle $\arc(\pi)$, and define the vector of
phase differences $\delta \triangleq B_c^{T} \theta =
(\theta_{2}-\theta_{1},\dots) \in {[-\pi,+\pi]}^{n(n-1)/2}$. By taking the
derivative $d/dt\, \delta(t)$ the phase differences satisfy
\begin{align}
  \dot \delta
  &=
  B_{c}^{T} \omega - B_{c}^{T} B \diag(\{a_{ij} \}_{\{i,j\} \in \mc E}) \sinbf(B^{T} \theta)
  \nonumber\\
  &=
  B_{c}^{T} \omega - B_{c}^{T} B_{c} \diag(\{a_{ij} \}_{i,j \in \until n, i < j}) \sinbf(\delta),
  \label{eq: delta dynamics}
\end{align}
where $\sinbf(x) = (\sin(x_{1}),\dots,\sin(x_{n}))$ for a vector $x \in
\real^{n}$. Notice that for $\theta(0) \in \arc(\pi)$ the
$\delta$-dynamics \eqref{eq: delta dynamics} are well-defined for an open
interval of time. In the following, we will show that the set $\{ \delta
\in \mathbb R^{n}:\, \| \delta \|_{2} < \subscr{\gamma}{max} \}$ is
positively invariant under condition \eqref{eq: frequ sync condition -
  1}. As a consequence, the set $\{ \delta \in \mathbb R^{n}:\, \| \delta
\|_{\infty} < \subscr{\gamma}{max} \leq \pi \}$ is positively invariant as
well, and the $\delta$-coordinates are well defined for all $t \geq 0$.

The Lyapunov function \eqref{eq: two-norm Lyapunov function} reads in
$\delta$-coordinates as $W(\delta) = \frac{1}{2} \| \delta \|^{2}$, and its
derivative along trajectories of \eqref{eq: delta dynamics}\,is
\begin{align}
  \dot W(\delta) &= \delta^{T} B_{c}^{T} \omega - \delta^{T}B_{c}^{T} B_{c}
  \diag(\{a_{ij} \}_{i < j}) \sinbf(\delta) \nonumber\\
  &= \delta^{T} B_{c}^{T} \omega - n\,\delta^{T} \diag(\{a_{ij} \}_{i < j})
  \sinbf(\delta)
  \label{eq: derivative of two-norm Lyapunov function - 1}
  \,,
\end{align}
where the second equality follows from the identity
\begin{equation*}
	\delta^{T}B_{c}^{T} B_{c}
	\!=\!
	\theta^{T}B_{c}B_{c}^{T} B_{c}
	\!=\!
	\theta^{T} (n I_{n} -  \fvec 1_{n\times n}) B_{c}
	\!=\! n\theta^{T} B_{c}
	\!=\! n\delta.
\end{equation*}
For $\|\delta_{2}\| \leq \rho$, $\rho \in {[0,\pi[}$, consider the following inequalities
\begin{multline*}
	n\,\delta^{T} \diag(\{a_{ij} \}_{i < j}) \sinbf(\delta)\\
	=
	n\,(B_{c}^{T}\theta)^{T} \diag(\{a_{ij} \sinc(\theta_{i}-\theta_{j}) \}_{i < j}) (B_{c}^{T}\theta)\\
	\geq 
	n\sinc(\rho) \, (B_{c}^{T}\theta)^{T} \diag(\{a_{ij} \}_{i < j}) (B_{c}^{T}\theta)\\
	\geq
	\lambda_{2}(L) \sinc(\rho) \| B_{c}^{T} \theta \|_{2}^{2}
	=
	\lambda_{2}(L) \sinc(\rho) \| \delta \|_{2}^{2}
	\,,
\end{multline*}
where the last inequality follows from \cite[Lemma 4.7]{FD-FB:09z}. Hence, the derivative \eqref{eq: derivative of two-norm Lyapunov function - 1} simplifies further to
\begin{equation}
	\dot W(\delta)
	\leq
	\delta^{T} B_{c}^{T} \omega
	- \lambda_{2}(L) \sinc(\rho) \| \delta \|_{2}^{2}
	\label{eq: derivative of two-norm Lyapunov function - 2}
	.
\end{equation}
In the following we regard $B_{c}^{T}\omega$ as external disturbance
affecting the otherwise stable $\delta$-dynamics \eqref{eq: delta dynamics}
and apply ultimate boundedness arguments \cite{HKK:02}. Note that the
right-hand side of \eqref{eq: derivative of two-norm Lyapunov function - 2}
is strictly negative for
\begin{equation*}
  \|\delta \|_{2} > 
  \mu_{c} \triangleq 
  \frac{\|B_{c}^{T} \omega\|_{2}}{\lambda_{2}(L) \sinc(\rho)}
  =
  \frac{\subscr{\lambda}{critical}}{\lambda_{2}(L) \sinc(\rho)}
  \,.
\end{equation*}
Pick $\epsilon \in ]0,1[$. If $\rho \geq \|\delta \|_{2} \geq \mu_{c}/\epsilon$, then the right-hand side of \eqref{eq: derivative of two-norm Lyapunov function - 2} is upper-bounded by
\begin{equation*}
  \dot{W}(\delta) 
  \leq - \bigl(1-\epsilon) \cdot  \lambda_{2}(L) \sinc(\rho) W(\delta) \,. 
\end{equation*}
In the following, choose $\mu$ such that $\rho > \mu > \mu_{c}$ and let $\epsilon = \mu_{c}/\mu \in {]0,1[}$.
By standard ultimate boundedness arguments \cite[Theorem 4.18]{HKK:02}, for $\|\delta(0)\|_{2} \leq \rho$, there is $T \geq 0$ such that $\|\delta(t)\|_{2}$ is exponentially decaying for $t \in [0,T]$ and $\|\delta(t)\|_{2} \leq \mu $ for all $t\geq T$. 
For the choice $\mu = \gamma$ with $\gamma \in {[0,\pi/2[}$, the condition $\mu > \mu_{c}$ reduces to
\begin{equation}
 	\gamma\sinc(\rho)
  	>
	\subscr{\lambda}{critical}/\lambda_{2}(L)
  	\label{eq: mu > muc condition}
	\,.
\end{equation}
Now, we perform a final analysis of the bound \eqref{eq: mu > muc condition}. 
The left-hand side of \eqref{eq: mu > muc condition} is an increasing function of $\gamma$ and a decreasing function of $\rho$. Therefore, there exists some $(\rho,\gamma)$ in the convex set 
$\Lambda \triangleq \{ (\rho,\gamma) :\, \rho \in {[0,\pi[} \,,\;  \gamma \in {[0,\pi/2[} \,,\; \rho > \gamma \}$ 
satisfying equation \eqref{eq: mu > muc condition} if and only if the inequality \eqref{eq: mu > muc condition} is true at $\rho = \gamma = \pi/2$, where the left-hand side of \eqref{eq: mu > muc condition} achieves its supremum in $\Lambda$. The latter condition is equivalent to inequality \eqref{eq: frequ sync condition - 1}.  Additionally, if these two equivalent statements are true, then there is an open set of points in $\Lambda$ satisfying \eqref{eq: mu > muc condition}, which is bounded by the unique curve that satisfies inequality \eqref{eq: mu > muc condition} with the equality sign, namely $f(\rho,\gamma) = 0$, where $f: \Lambda \to \real$,
$f(\rho,\gamma) = \gamma \sinc(\rho) -  \subscr{\lambda}{critical}/\lambda_{2}(L)$. 
Consequently, for every $(\rho,\gamma) \in \{ (\rho,\gamma) \in \Lambda:\, f(\rho,\gamma) > 0 \}$, it follows for $\|\delta(0)\|_{2} \leq \rho$ that there is $T \geq 0$ such that $\|\delta(t)\|_{2} \leq \gamma$ for all $t \geq T$. The supremum value for $\rho$ is given by $\subscr{\rho}{max} \in {]\pi/2,\pi]}$ solving the equation $f(\subscr{\rho}{max},\pi/2) =0$ and the infimum value of $\gamma$ by $\subscr{\gamma}{min} \in \!{[0,\pi/2[}$ solving the equation $f(\subscr{\gamma}{min},\subscr{\gamma}{min})\!=\!0$. 

This proves statement 1) (where we replaced $\subscr{\rho}{max}$ by $\subscr{\gamma}{max}$) and shows that there is $T \geq 0$ such that $\|B_{c}^{T}\theta(t)\|_{\infty} \!\leq\! \|B_{c}^{T}\theta(t)\|_{2} \leq \subscr{\gamma}{min} < \pi/2$ for all $t \geq T$. Thus, $\theta(t) \in \bar\Delta_{G}(\subscr{\gamma}{min})$ for $t \geq T$, and frequency synchronization can be established analogously to the proof of Theorem \ref{Theorem: Kuramoto frequency sync}.

\subsection{Proof of Theorem \ref{Theorem: Synchronization in coupled oscillators II}}
\label{Subsection: Proof of freq sync in sparse graphs II}

According to Lemma \ref{Lemma: stable equilibrium in pi/2 arc}, there exists a locally exponentially stable synchronization manifold $[\theta] \in \bar\Delta_{G}(\gamma)$, $\gamma \in {[0,\pi/2[}$, if and only if there is an equilibrium $\theta \in \bar\Delta_{G}(\gamma)$. The equilibrium equations \eqref{eq: equilibrium equations - component-wise} can be rewritten\,as 
\begin{equation}
	\omega = L(B^{T}\theta) \theta
	\,,
	\label{eq: equilibrium equations - vector-wise}
\end{equation}
where $L(B^{T}\theta) = B \diag(\{a_{ij}\sinc(\theta_{i} - \theta_{j})\}_{\{i,j\} \in \mc E})B^{T}$ is the Laplacian matrix associated with the graph $G(\mc V,\mc E,\tilde A)$  with nonnegative edge weights $\tilde a_{ij} = a_{ij}\sinc(\theta_{i} - \theta_{j})$ for $\theta \in \bar\Delta_{G}(\gamma)$. Since for any weighted Laplacian matrix $L$, we have that $L \cdot L^{\dagger} = L^{\dagger} \cdot L = I_{n} - (1/n) \fvec 1_{n \times n}$ (follows from the singular value decomposition \cite{FD-FB:11d}), a multiplication of  equation \eqref{eq: equilibrium equations - vector-wise} from the left by $B^{T} L(B^{T}\theta)^{\dagger}$ yields
\begin{equation}
	B^{T} L(B^{T}\theta)^{\dagger} \omega = B^{T}\theta
	\label{eq: fixed-point equations}
	\,.
\end{equation}
Note that the left-hand side of equation \eqref{eq: fixed-point equations} is a\,continuous%
\footnote{
The continuity can be established when re-writing equations \eqref{eq: equilibrium equations - vector-wise} and \eqref{eq: fixed-point equations} in the quotient space $\fvec 1_{n}^{\perp}$, where $L(B^{T}\theta)$ is nonsingular, and using the fact that the inverse of a matrix is a continuous function of its elements.
}
function for  
$\theta \in \bar\Delta_{G}(\gamma)$. 
Consider the formal substitution $x = B^{T} \theta$, the compact and convex
set $\mc S_{\infty}(\gamma) = \{x \in B^{T}\mathbb R^{n}:\, \|x\|_{\infty}
\leq \gamma \}$, and the continuous map\,$f: \mc S_{\infty}(\gamma) \to
\real$ given by $f(x) = B^{T} L(x)^{\dagger} \omega$. Then equation
\eqref{eq: fixed-point equations} reads as the fixed-point equation $f(x) =
x$, and we can invoke {\em Brouwers's Fixed Point Theorem} which states
that every continuous map from a compact and convex set to itself has a
fixed point, see for instance \cite[Section 7, Corollary 8]{EHS:94}.

Since the analysis of the map $f$ in the $\infty$-norm is very hard in the
general case, we resort to a $2$-norm analysis and restrict ourselves to
the set $\mc S_{2}(\gamma) = \{x \in B^{T}\mathbb R^{n}:\, \|x\|_{2} \leq
\gamma \} \subseteq \mc S_{\infty}(\gamma)$.  By Brouwer's Fixed Point
Theorem, there exists a solution $x \in \mc S_{2}(\gamma)$ to the equations
$x = f(x)$ if and only if $\|f(x)\|_{2} \leq \gamma$ for all $x \in \mc
S_{2}(\gamma)$, or equivalently if and only if
\begin{equation}
	\max_{x\in \mc S_{2}(\gamma)} \left\| B^{T} L(x)^{\dagger} \omega \right\|_{2} \leq \gamma
	\,.
	\label{eq: Brouwer condition}
\end{equation}
In the following we show that \eqref{eq: frequ sync condition - 2} is a
sufficient condition for inequality \eqref{eq: Brouwer condition}. 
%

First, we establish some identities.
For a Laplacian matrix $L$, we obtain $L^{\dagger} = V \diag(0,\{1/\lambda_{i}(L) \}_{i=2,\dots,n}) V^{T}$, where $\lambda_{1}(L) = 0$ and $\lambda_{i}(L) > 0$, $i \in \{2,\dots,n\}$, are the eigenvalues of $L$ and $V \in \real^{n \times n}$ is an associated orthonormal matrix of eigenvectors. It follows that $V\diag\left( 0,1,\dots,1 \right) V^{T} = I_{n} - (1/n) \fvec 1_{n \times n}$, and since $\omega \perp \fvec 1_{n}$, there exists $\alpha \in \real^{|\mc E|}$ (not necessarily unique), such that $\omega = B\alpha$. By means of these identities, the left-hand side of \eqref{eq: Brouwer condition} can be simplified and upper-bounded for all $x\in \mc S_{2}(\gamma)$:
\begin{align}
	&
	\left\| B^{T} L(x)^{\dagger} \omega \right\|_{2}
	=
	\left\| B^{T} L(x)^{\dagger}  B \alpha \right\|_{2}
	=\nonumber\\
	&
	\left\| B^{T} V(x) \diag\left( 0,\frac{1}{\lambda_{2}(L(x))},\dots, \frac{1}{\lambda_{n}(L(x))}  \right) V^{T}(x)  B \alpha \right\|_{2}
	\nonumber\\
	&\leq
	\frac{1}{\lambda_{2}(L(x))} \cdot \left\| B^{T} V(x) \diag\left( 0,1,\dots, 1  \right) V^{T}(x)  B \alpha \right\|_{2}
	\nonumber\\
	&\quad=
	(1/\lambda_{2}(L(x))) \cdot \left\| B^{T} \omega \right\|_{2} 
	\label{eq: bounding id}
	\,.
\end{align}
Thus, a sufficient condition for inequality \eqref{eq: Brouwer condition} to be true can be derived as follows:
\begin{align*}
	&
	\max_{x\in \mc S_{2}(\gamma)} \left\| B^{T} L(x)^{\dagger} \omega \right\|_{2}
	\leq
	\left\| B^{T} \omega \right\|_{2} \max_{x\in \mc S_{2}(\gamma)} \left(1/\lambda_{2}\bigl(L(x)\bigr))\right)
	\\
	&\quad\leq
	\left\| B^{T}\omega \right\|_{2} \max_{x \in \{x \in \mathbb R^{|\mc E|}:\, \|x\|_{\infty} \leq \gamma \}} \left( 1/ \lambda_{2}\bigl(L(x)\bigr)) \right)
	\\&\quad= 
	\left\| B^{T}\omega \right\|_{2} / \left( \lambda_{2}(L) \cdot \sinc(\gamma) \right)
	\overset{!}{\leq} \gamma
	\,,
\end{align*}
where we used identity \eqref{eq: bounding id}, we enlarged the domain $\mc S_{2}(\gamma)$ to $\{x \in
\mathbb R^{|\mc E|}:\, \|x\|_{\infty} \leq \gamma \}$, and we used the fact
$\lambda_{2}(L(x)) \geq \lambda_{2}(L) \cdot \sinc(\gamma)$ for
$\|x\|_{\infty} \leq \gamma$.
In summary,\,we conclude that there is a locally exponentially stable synchronization manifold 
$[\theta] \in \{ \theta \in \torus^{n}: \| B^{T} \theta\|_{2} \!\leq\! \gamma\} \subseteq \bar\Delta_{G}(\gamma)$\,if
\begin{equation}
	\lambda_{2}(L) \sin(\gamma) \geq \| B^{T} \omega \|_{2}
	\label{eq: frequ sync condition - 2 - gamma}
	\,.
\end{equation}
Since the left-hand side of \eqref{eq: frequ sync condition - 2 - gamma} is a concave function of $\gamma \in {[0,\pi/2[}$, there exists an open set of $\gamma \in {[0,\pi/2[}$\,satisfying equation \eqref{eq: frequ sync condition - 2 - gamma} if and only if equation \eqref{eq: frequ sync condition - 2 - gamma} is true with the strict equality sign at $\gamma^{*}=\pi/2$, which corresponds to condition \eqref{eq: frequ sync condition - 2}. Additionally, if these two equivalent  statements are true, then there exists a unique $\subscr{\gamma}{min}\in {[0,\pi/2[}$ that satisfies equation \eqref{eq: critical bound} with the equality sign, namely $\sin(\subscr{\gamma}{min}) = \| B^{T}\omega \|_{2} / \lambda_{2}(L)$. This concludes the proof.

\end{appendix}


\bibliographystyle{IEEEtran}
\bibliography{alias,Main,FB}

\begin{thebibliography}{100}
\providecommand{\url}[1]{#1}
\csname url@samestyle\endcsname
\providecommand{\newblock}{\relax}
\providecommand{\bibinfo}[2]{#2}
\providecommand{\BIBentrySTDinterwordspacing}{\spaceskip=0pt\relax}
\providecommand{\BIBentryALTinterwordstretchfactor}{4}
\providecommand{\BIBentryALTinterwordspacing}{\spaceskip=\fontdimen2\font plus
\BIBentryALTinterwordstretchfactor\fontdimen3\font minus
  \fontdimen4\font\relax}
\providecommand{\BIBforeignlanguage}[2]{{%
\expandafter\ifx\csname l@#1\endcsname\relax
\typeout{** WARNING: IEEEtran.bst: No hyphenation pattern has been}%
\typeout{** loaded for the language `#1'. Using the pattern for}%
\typeout{** the default language instead.}%
\else
\language=\csname l@#1\endcsname
\fi
#2}}
\providecommand{\BIBdecl}{\relax}
\BIBdecl

\bibitem{CH:1673}
C.~Huygens, \emph{Horologium Oscillatorium}, Paris, France, 1673.

\bibitem{SHS:03}
S.~H. Strogatz, \emph{SYNC: The Emerging Science of Spontaneous Order}.\hskip
  1em plus 0.5em minus 0.4em\relax Hyperion, 2003.

\bibitem{ATW:01}
A.~T. Winfree, \emph{The Geometry of Biological Time}, 2nd~ed.\hskip 1em plus
  0.5em minus 0.4em\relax Springer, 2001.

\bibitem{ATW:67}
------, ``Biological rhythms and the behavior of populations of coupled
  oscillators,'' \emph{Journal of Theoretical Biology}, vol.~16, no.~1, pp.
  15--42, 1967.

\bibitem{YK:75}
Y.~Kuramoto, ``Self-entrainment of a population of coupled non-linear
  oscillators,'' in \emph{Int. Symposium on Mathematical Problems in
  Theoretical Physics}, ser. Lecture Notes in Physics, H.~Araki, Ed.\hskip 1em
  plus 0.5em minus 0.4em\relax Springer, 1975, vol.~39, pp. 420--422.

\bibitem{YK:84}
------, \emph{Chemical Oscillations, Waves, and Turbulence}.\hskip 1em plus
  0.5em minus 0.4em\relax Springer, 1984.

\bibitem{SHS:00}
S.~H. Strogatz, ``From {K}uramoto to {C}rawford: Exploring the onset of
  synchronization in populations of coupled oscillators,'' \emph{Physica D:
  Nonlinear Phenomena}, vol. 143, no.~1, pp. 1--20, 2000.

\bibitem{JAA-LLB-CJPV-FR-RS:05}
J.~A. Acebr{\'o}n, L.~L. Bonilla, C.~J.~P. Vicente, F.~Ritort, and R.~Spigler,
  ``The {K}uramoto model: {A} simple paradigm for synchronization phenomena,''
  \emph{Reviews of Modern Physics}, vol.~77, no.~1, pp. 137--185, 2005.

\bibitem{DCM-EPM-JJ:87}
D.~C. Michaels, E.~P. Matyas, and J.~Jalife, ``Mechanisms of sinoatrial
  pacemaker synchronization: a new hypothesis,'' \emph{Circulation Research},
  vol.~61, no.~5, pp. 704--714, 1987.

\bibitem{CL-DRW-SHS-RSM:97}
C.~Liu, D.~R. Weaver, S.~H. Strogatz, and S.~M. Reppert, ``Cellular
  construction of a circadian clock: period determination in the
  suprachiasmatic nuclei,'' \emph{Cell}, vol.~91, no.~6, pp. 855--860, 1997.

\bibitem{FV-JPL-ER-JM:01}
F.~Varela, J.~P. Lachaux, E.~Rodriguez, and J.~Martinerie, ``The brainweb:
  {P}hase synchronization and large-scale integration,'' \emph{Nature Reviews
  Neuroscience}, vol.~2, no.~4, pp. 229--239, 2001.

\bibitem{EB-PH-JM:03}
E.~Brown, P.~Holmes, and J.~Moehlis, ``Globally coupled oscillator networks,''
  in \emph{Perspectives and Problems in Nonlinear Science: A Celebratory Volume
  in Honor of Larry Sirovich}, E.~Kaplan, J.~E. Marsden, and K.~R. Sreenivasan,
  Eds.\hskip 1em plus 0.5em minus 0.4em\relax Springer, 2003, pp. 183--215.

\bibitem{SMC-GBE-MCV-JMB:97}
S.~M. Crook, G.~B. Ermentrout, M.~C. Vanier, and J.~M. Bower, ``The role of
  axonal delay in the synchronization of networks of coupled cortical
  oscillators,'' \emph{Journal of Computational Neuroscience}, vol.~4, no.~2,
  pp. 161--172, 1997.

\bibitem{AKG-BC-EKP:71}
A.~K. Ghosh, B.~Chance, and E.~K. Pye, ``Metabolic coupling and synchronization
  of {NADH} oscillations in yeast cell populations,'' \emph{Archives of
  Biochemistry and Biophysics}, vol. 145, no.~1, pp. 319--331, 1971.

\bibitem{JB:88}
J.~Buck, ``Synchronous rhythmic flashing of fireflies. {II}.'' \emph{Quarterly
  Review of Biology}, vol.~63, no.~3, pp. 265--289, 1988.

\bibitem{TJW:69}
T.~J. Walker, ``Acoustic synchrony: two mechanisms in the snowy tree cricket,''
  \emph{Science}, vol. 166, no. 3907, pp. 891--894, 1969.

\bibitem{NK-GBE:88}
N.~Kopell and G.~B. Ermentrout, ``Coupled oscillators and the design of central
  pattern generators,'' \emph{Mathematical Biosciences}, vol.~90, no. 1-2, pp.
  87--109, 1988.

\bibitem{NEL-TS-NB-LS-IDC-SAL:12}
N.~E. Leonard, T.~Shen, B.~Nabet, L.~Scardovi, I.~D. Couzin, and S.~A. Levin,
  ``Decision versus compromise for animal groups in motion,'' \emph{Proceedings
  of the National Academy of Sciences}, vol. 109, no.~1, pp. 227--232, 2012.

\bibitem{DAP-NEL-RS-DG-JKP:07}
D.~A. Paley, N.~E. Leonard, R.~Sepulchre, D.~Grunbaum, and J.~K. Parrish,
  ``Oscillator models and collective motion,'' \emph{{IEEE} Control Systems
  Magazine}, vol.~27, no.~4, pp. 89--105, 2007.

\bibitem{ZN-ER-TV-YB-AIB:00}
Z.~N{\'e}da, E.~Ravasz, T.~Vicsek, Y.~Brechet, and A.~L. Barab{\'a}si,
  ``Physics of the rhythmic applause,'' \emph{Physical Review E}, vol.~61,
  no.~6, p. 6987, 2000.

\bibitem{HD:92}
H.~Daido, ``Quasientrainment and slow relaxation in a population of oscillators
  with random and frustrated interactions,'' \emph{Physical Review Letters},
  vol.~68, no.~7, pp. 1073--1076, 1992.

\bibitem{GJ-JA-DB-ACCC-CPV:01}
G.~Jongen, J.~Anem{\"u}ller, D.~Boll{\'e}, A.~C.~C. Coolen, and
  C.~Perez-Vicente, ``Coupled dynamics of fast spins and slow exchange
  interactions in the {XY} spin glass,'' \emph{Journal of Physics A:
  Mathematical and General}, vol.~34, no.~19, pp. 3957--3984, 2001.

\bibitem{JP:98}
J.~Pantaleone, ``Stability of incoherence in an isotropic gas of oscillating
  neutrinos,'' \emph{Physical Review D}, vol.~58, no.~7, p. 073002, 1998.

\bibitem{KW-PC-SHS:98}
K.~Wiesenfeld, P.~Colet, and S.~H. Strogatz, ``Frequency locking in {J}osephson
  arrays: {C}onnection with the {K}uramoto model,'' \emph{Physical Review E},
  vol.~57, no.~2, pp. 1563--1569, 1998.

\bibitem{IZK-YZ-JLH:02}
I.~Z. Kiss, Y.~Zhai, and J.~L. Hudson, ``Emerging coherence in a population of
  chemical oscillators,'' \emph{Science}, vol. 296, no. 5573, p. 1676, 2002.

\bibitem{PAT:03}
P.~A. Tass, ``A model of desynchronizing deep brain stimulation with a
  demand-controlled coordinated reset of neural subpopulations,''
  \emph{Biological Cybernetics}, vol.~89, no.~2, pp. 81--88, 2003.

\bibitem{AN-JM:11}
A.~Nabi and J.~Moehlis, ``Single input optimal control for globally coupled
  neuron networks,'' \emph{Journal of Neural Engineering}, vol.~8, p. 065008,
  2011.

\bibitem{RS-DP-NEL:07}
R.~Sepulchre, D.~A. Paley, and N.~E. Leonard, ``Stabilization of planar
  collective motion: {A}ll-to-all communication,'' \emph{IEEE Transactions on
  Automatic Control}, vol.~52, no.~5, pp. 811--824, 2007.

\bibitem{RS-DP-NEL:08}
------, ``Stabilization of planar collective motion with limited
  communication,'' \emph{IEEE Transactions on Automatic Control}, vol.~53,
  no.~3, pp. 706--719, 2008.

\bibitem{DJK:08}
D.~J. Klein, ``Coordinated control and estimation for multi-agent systems:
  Theory and practice,'' Ph.D. dissertation, University of Washington, 2008.

\bibitem{DJK-PL-KAM-TJ:08}
D.~J. Klein, P.~Lee, K.~A. Morgansen, and T.~Javidi, ``Integration of
  communication and control using discrete time {K}uramoto models for
  multivehicle coordination over broadcast networks,'' \emph{IEEE Journal on
  Selected Areas in Communications}, vol.~26, no.~4, pp. 695--705, 2008.

\bibitem{RMMU-RM-SD:11}
M.~M.~U. Rahman, R.~Mudumbai, and S.~Dasgupta, ``Consensus based carrier
  synchronization in a two node network,'' in \emph{{IFAC} {W}orld {C}ongress},
  Milan, Italy, Aug. 2011, pp. 10\,038--10\,043.

\bibitem{GK-AGV-PM:00}
G.~Kozyreff, A.~G. Vladimirov, and P.~Mandel, ``Global coupling with time delay
  in an array of semiconductor lasers,'' \emph{Physical Review Letters},
  vol.~85, no.~18, pp. 3809--3812, 2000.

\bibitem{FCH-EMI:00}
F.~C. Hoppensteadt and E.~M. Izhikevich, ``Synchronization of laser
  oscillators, associative memory, and optical neurocomputing,'' \emph{Physical
  Review E}, vol.~62, no.~3, pp. 4010--4013, 2000.

\bibitem{RAY-RCC:02}
R.~A. York and R.~C. Compton, ``Quasi-optical power combining using mutually
  synchronized oscillator arrays,'' \emph{IEEE Transactions on Microwave Theory
  and Techniques}, vol.~39, no.~6, pp. 1000--1009, 2002.

\bibitem{WCL-FG-WCH-KD:85}
W.~C. Lindsey, F.~Ghazvinian, W.~C. Hagmann, and K.~Dessouky, ``Network
  synchronization,'' \emph{Proceedings of the IEEE}, vol.~73, no.~10, pp.
  1445--1467, 1985.

\bibitem{OS-US-YBN-SS:08}
O.~Simeone, U.~Spagnolini, Y.~Bar-Ness, and S.~H. Strogatz, ``Distributed
  synchronization in wireless networks,'' \emph{IEEE Signal Processing
  Magazine}, vol.~25, no.~5, pp. 81--97, 2008.

\bibitem{YWH-AS:05}
Y.~W. Hong and A.~Scaglione, ``A scalable synchronization protocol for large
  scale sensor networks and its applications,'' \emph{IEEE Journal on Selected
  Areas in Communications}, vol.~23, no.~5, pp. 1085--1099, 2005.

\bibitem{RB-AC-LQ-SS-STP:10}
R.~Baldoni, A.~Corsaro, L.~Querzoni, S.~Scipioni, and S.~T. Piergiovanni,
  ``Coupling-based internal clock synchronization for large-scale dynamic
  distributed systems,'' \emph{IEEE Transactions on Parallel and Distributed
  Systems}, vol.~21, no.~5, pp. 607--619, 2010.

\bibitem{YW-NF-JFD:12}
Y.~Wang, F.~N{\'u}{\~n}ez, and F.~J. Doyle, ``Increasing sync rate of
  pulse-coupled oscillators via phase response function design: theory and
  application to wireless networks,'' \emph{IEEE Transactions on Control
  Systems Technology}, 2012, to appear.

\bibitem{EM-AT:11b}
E.~Mallada and A.~Tang, ``Distributed clock synchronization: Joint frequency
  and phase consensus,'' in \emph{{IEEE} Conf.\ on Decision and Control and
  European Control Conference}, Orlando, FL, USA, Dec. 2011, pp. 6742--6747.

\bibitem{SB-SG:07}
S.~Barbarossa and G.~Scutari, ``Decentralized maximum-likelihood estimation for
  sensor networks composed of nonlinearly coupled dynamical systems,''
  \emph{IEEE Transactions on Signal Processing}, vol.~55, no.~7, pp.
  3456--3470, 2007.

\bibitem{JWSP-FD-FB:12j}
J.~W. Simpson-Porco, F.~D{\"o}rfler, and F.~Bullo, ``Droop-controlled inverters
  are {Kuramoto} oscillators,'' in \emph{IFAC Workshop on Distributed
  Estimation and Control in Networked Systems}, Santa Barbara, CA, USA, Sep.
  2012, to appear.

\bibitem{AA-ADG-JK-YM-CZ:08}
A.~Arenas, A.~D{\'\i}az-Guilera, J.~Kurths, Y.~Moreno, and C.~Zhou,
  ``Synchronization in complex networks,'' \emph{Physics Reports}, vol. 469,
  no.~3, pp. 93--153, 2008.

\bibitem{SB-VL-YM-MC-DUH:06}
S.~Boccaletti, V.~Latora, Y.~Moreno, M.~Chavez, and D.~U. Hwang, ``Complex
  networks: {S}tructure and dynamics,'' \emph{Physics Reports}, vol. 424, no.
  4-5, pp. 175--308, 2006.

\bibitem{SHS:01}
S.~H. Strogatz, ``Exploring complex networks,'' \emph{Nature}, vol. 410, no.
  6825, pp. 268--276, 2001.

\bibitem{JAKS-GVO:08}
J.~A.~K. Suykens and G.~V. Osipov, ``Introduction to focus issue:
  {S}ynchronization in complex networks,'' \emph{Chaos}, vol.~18, no.~3, pp.
  037\,101--037\,101, 2008.

\bibitem{ROS-JAF-RMM:07}
R.~Olfati-Saber, J.~A. Fax, and R.~M. Murray, ``Consensus and cooperation in
  networked multi-agent systems,'' \emph{Proceedings of the IEEE}, vol.~95,
  no.~1, pp. 215--233, 2007.

\bibitem{WR-RWB-EMA:07}
W.~Ren, R.~W. Beard, and E.~M. Atkins, ``Information consensus in multivehicle
  cooperative control: {C}ollective group behavior through local interaction,''
  \emph{{IEEE} Control Systems Magazine}, vol.~27, no.~2, pp. 71--82, 2007.

\bibitem{FB-JC-SM:09}
F.~Bullo, J.~Cort{\'e}s, and S.~Mart{\'\i}nez, \emph{Distributed Control of
  Robotic Networks}, ser. Applied Mathematics Series.\hskip 1em plus 0.5em
  minus 0.4em\relax Princeton University Press, 2009.

\bibitem{LM:05}
L.~Moreau, ``Stability of multiagent systems with time-dependent communication
  links,'' \emph{IEEE Transactions on Automatic Control}, vol.~50, no.~2, pp.
  169--182, 2005.

\bibitem{AJ-NM-MB:04}
A.~Jadbabaie, N.~Motee, and M.~Barahona, ``On the stability of the {K}uramoto
  model of coupled nonlinear oscillators,'' in \emph{{A}merican {C}ontrol
  {C}onference}, Boston, MA, USA, Jun. 2004, pp. 4296--4301.

\bibitem{NC-MWS:08}
N.~Chopra and M.~W. Spong, ``On exponential synchronization of {K}uramoto
  oscillators,'' \emph{IEEE Transactions on Automatic Control}, vol.~54, no.~2,
  pp. 353--357, 2009.

\bibitem{ZL-BF-MM:07}
Z.~Lin, B.~Francis, and M.~Maggiore, ``State agreement for continuous-time
  coupled nonlinear systems,'' \emph{SIAM Journal on Control and Optimization},
  vol.~46, no.~1, pp. 288--307, 2007.

\bibitem{AS-RS:07}
A.~Sarlette and R.~Sepulchre, ``Consensus optimization on manifolds,''
  \emph{SIAM Journal on Control and Optimization}, vol.~48, no.~1, pp. 56--76,
  2009.

\bibitem{LS-AS-RS:06}
L.~Scardovi, A.~Sarlette, and R.~Sepulchre, ``Synchronization and balancing on
  the {$N$}-torus,'' \emph{Systems \& Control Letters}, vol.~56, no.~5, pp.
  335--341, 2007.

\bibitem{ROS:06b}
R.~Olfati-Saber, ``Swarms on sphere: {A} programmable swarm with synchronous
  behaviors like oscillator networks,'' in \emph{{IEEE} Conf.\ on Decision and
  Control}, San Diego, CA, USA, 2006, pp. 5060--5066.

\bibitem{GBE:91}
G.~B. Ermentrout, ``An adaptive model for synchrony in the firefly
  \emph{pteroptyx malaccae},'' \emph{Journal of Mathematical Biology}, vol.~29,
  no.~6, pp. 571--585, 1991.

\bibitem{SYH-EJ-MJK:10}
S.~Y. Ha, E.~Jeong, and M.~J. Kang, ``Emergent behaviour of a generalized
  {V}iscek-type flocking model,'' \emph{Nonlinearity}, vol.~23, no.~12, pp.
  3139--3156, 2010.

\bibitem{SYH-CL-BR-MS:11}
S.~Ha, C.~Lattanzio, B.~Rubino, and M.~Slemrod, ``Flocking and synchronization
  of particle models,'' \emph{Quarterly Applied Mathematics}, vol.~69, pp.
  91--103, 2011.

\bibitem{ARB-DJH:81}
A.~R. Bergen and D.~J. Hill, ``A structure preserving model for power system
  stability analysis,'' \emph{IEEE Transactions on Power Apparatus and
  Systems}, vol. 100, no.~1, pp. 25--35, 1981.

\bibitem{PWS-MAP:98}
P.~W. Sauer and M.~A. Pai, \emph{Power System Dynamics and Stability}.\hskip
  1em plus 0.5em minus 0.4em\relax Prentice Hall, 1998.

\bibitem{HDC-CCC-GC:95}
H.-D. Chiang, C.~C. Chu, and G.~Cauley, ``Direct stability analysis of electric
  power systems using energy functions: {T}heory, applications, and
  perspective,'' \emph{Proceedings of the IEEE}, vol.~83, no.~11, pp.
  1497--1529, 1995.

\bibitem{FD-FB:09z}
F.~D{\"o}rfler and F.~Bullo, ``Synchronization and transient stability in power
  networks and non-uniform {K}uramoto oscillators,'' \emph{SIAM Journal on
  Control and Optimization}, vol.~50, no.~3, pp. 1616--1642, 2012.

\bibitem{JP:02}
J.~Pantaleone, ``Synchronization of metronomes,'' \emph{American Journal of
  Physics}, vol.~70, p. 992, 2002.

\bibitem{SHS-DMA-AMR-BE-EO:05}
S.~H. Strogatz, D.~M. Abrams, A.~McRobie, B.~Eckhardt, and E.~Ott,
  ``Theoretical mechanics: {C}rowd synchrony on the {Millennium Bridge},''
  \emph{Nature}, vol. 438, no. 7064, pp. 43--44, 2005.

\bibitem{MB-MFS-HR-KW:02}
M.~Bennett, M.~F. Schatz, H.~Rockwood, and K.~Wiesenfeld, ``Huygens's clocks,''
  \emph{Proceedings: Mathematical, Physical and Engineering Sciences}, vol.
  458, no. 2019, pp. 563--579, 2002.

\bibitem{YPC-SYH-SBH:11}
Y.-P. Choi, S.-Y. Ha, and S.-B. Yun, ``Complete synchronization of {K}uramoto
  oscillators with finite inertia,'' \emph{Physica D}, vol. 240, no.~1, pp.
  32--44, 2011.

\bibitem{HAT-AJL-SO:97}
H.~A. Tanaka, A.~J. Lichtenberg, and S.~Oishi, ``Self-synchronization of
  coupled oscillators with hysteretic responses,'' \emph{Physica D: Nonlinear
  Phenomena}, vol. 100, no. 3-4, pp. 279--300, 1997.

\bibitem{HAT-AJL-SO:97b}
------, ``First order phase transition resulting from finite inertia in coupled
  oscillator systems,'' \emph{Physical Review Letters}, vol.~78, no.~11, pp.
  2104--2107, 1997.

\bibitem{HH-GSJ-MYC:02}
H.~Hong, G.~S. Jeon, and M.~Y. Choi, ``Spontaneous phase oscillation induced by
  inertia and time delay,'' \emph{Physical Review E}, vol.~65, no.~2, p.
  026208, 2002.

\bibitem{HH-MYC-JY-KSS:99}
H.~Hong, M.~Y. Choi, J.~Yi, and K.~S. Soh, ``Inertia effects on periodic
  synchronization in a system of coupled oscillators,'' \emph{Physical Review
  E}, vol.~59, no.~1, p. 353, 1999.

\bibitem{JAA-LLB-RS:00}
J.~A. Acebr{\'o}n, L.~L. Bonilla, and R.~Spigler, ``Synchronization in
  populations of globally coupled oscillators with inertial effects,''
  \emph{Physical Review E}, vol.~62, no.~3, pp. 3437--3454, 2000.

\bibitem{FD-FB:10w}
F.~D{\"o}rfler and F.~Bullo, ``On the critical coupling for {K}uramoto
  oscillators,'' \emph{SIAM Journal on Applied Dynamical Systems}, vol.~10,
  no.~3, pp. 1070--1099, 2011.

\bibitem{FDS-DA:07}
F.~{De~Smet} and D.~Aeyels, ``Partial entrainment in the finite
  {Kuramoto--Sakaguchi} model,'' \emph{Physica D: Nonlinear Phenomena}, vol.
  234, no.~2, pp. 81--89, 2007.

\bibitem{RM-SHS:07}
R.~Mirollo and S.~H. Strogatz, ``The spectrum of the partially locked state for
  the {K}uramoto model,'' \emph{Journal of Nonlinear Science}, vol.~17, no.~4,
  pp. 309--347, 2007.

\bibitem{YLM-OVP-PAT:05}
Y.~L. Maistrenko, O.~V. Popovych, and P.~A. Tass, ``Desynchronization and chaos
  in the {K}uramoto model,'' in \emph{Dynamics of Coupled Map Lattices and of
  Related Spatially Extended Systems}, ser. Lecture Notes in Physics, J.-R.
  Chazottes and B.~Fernandez, Eds.\hskip 1em plus 0.5em minus 0.4em\relax
  Springer, 2005, vol. 671, pp. 285--306.

\bibitem{RT:07}
R.~T{\"o}njes, ``Pattern formation through synchronization in systems of
  nonidentical autonomous oscillators,'' Ph.D. dissertation, Universit{\"a}ts
  Potsdam, Germany, 2007.

\bibitem{OVP-YLM-PAT:05}
O.~V. Popovych, Y.~L. Maistrenko, and P.~A. Tass, ``Phase chaos in coupled
  oscillators,'' \emph{Physical Review E}, vol.~71, no.~6, p. 065201, 2005.

\bibitem{JL:11}
J.~Lunze, ``Complete synchronization of {K}uramoto oscillators,'' \emph{Journal
  of Physics A: Mathematical and Theoretical}, vol.~44, p. 425102, 2011.

\bibitem{EC-PM:08}
E.~Canale and P.~Monz{\'o}n, ``Almost global synchronization of symmetric
  {K}uramoto coupled oscillators,'' in \emph{Systems Structure and
  Control}.\hskip 1em plus 0.5em minus 0.4em\relax InTech Education and
  Publishing, 2008, ch.~8, pp. 167--190.

\bibitem{MV-OM:09}
M.~Verwoerd and O.~Mason, ``On computing the critical coupling coefficient for
  the {K}uramoto model on a complete bipartite graph,'' \emph{SIAM Journal on
  Applied Dynamical Systems}, vol.~8, no.~1, pp. 417--453, 2009.

\bibitem{SYH-TH-KJH:10}
S.-Y. Ha, T.~Ha, and J.-H. Kim, ``On the complete synchronization of the
  {K}uramoto phase model,'' \emph{Physica D: Nonlinear Phenomena}, vol. 239,
  no.~17, pp. 1692--1700, 2010.

\bibitem{REM-SHS:05}
R.~E. Mirollo and S.~H. Strogatz, ``The spectrum of the locked state for the
  {K}uramoto model of coupled oscillators,'' \emph{Physica D: Nonlinear
  Phenomena}, vol. 205, no. 1-4, pp. 249--266, 2005.

\bibitem{DA-JAR:04}
D.~Aeyels and J.~A. Rogge, ``Existence of partial entrainment and stability of
  phase locking behavior of coupled oscillators,'' \emph{Progress on
  Theoretical Physics}, vol. 112, no.~6, pp. 921--942, 2004.

\bibitem{MV-OM:08}
M.~Verwoerd and O.~Mason, ``Global phase-locking in finite populations of
  phase-coupled oscillators,'' \emph{SIAM Journal on Applied Dynamical
  Systems}, vol.~7, no.~1, pp. 134--160, 2008.

\bibitem{GBE:85}
G.~B. Ermentrout, ``Synchronization in a pool of mutually coupled oscillators
  with random frequencies,'' \emph{Journal of Mathematical Biology}, vol.~22,
  no.~1, pp. 1--9, 1985.

\bibitem{VM-OM:11}
M.~Verwoerd and O.~Mason, ``A convergence result for the {K}uramoto model with
  all-to-all coupling,'' \emph{SIAM Journal on Applied Dynamical Systems},
  vol.~10, no.~3, pp. 906--920, 2011.

\bibitem{JCB-LDV-MJP:11}
J.~C. Bronski, L.~DeVille, and M.~J. Park, ``Fully synchronous solutions and
  the synchronization phase transition for the finite $n$ {K}uramoto model,''
  \emph{Arxiv preprint arXiv:1111.5302}, 2011.

\bibitem{JAR-DA:04}
J.~A. Rogge and D.~Aeyels, ``Stability of phase locking in a ring of
  unidirectionally coupled oscillators,'' \emph{Journal of Physics A}, vol.~37,
  pp. 11\,135--11\,148, 2004.

\bibitem{LDV:11}
L.~DeVille, ``Transitions amongst synchronous solutions for the stochastic
  {K}uramoto model,'' \emph{Nonlinearity}, vol.~25, no.~5, pp. 1--20, 2011.

\bibitem{UM-AP-FA:09}
U.~M{\"u}nz, A.~Papachristodoulou, and F.~Allg{\"o}wer, ``Consensus reaching in
  multi-agent packet-switched networks with non-linear coupling,''
  \emph{International Journal of Control}, vol.~82, no.~5, pp. 953--969, 2009.

\bibitem{EM-AT:10}
E.~Mallada and A.~Tang, ``Synchronization of phase-coupled oscillators with
  arbitrary topology,'' in \emph{{A}merican {C}ontrol {C}onference}, Baltimore,
  MD, USA, Jun. 2010, pp. 1777--1782.

\bibitem{LS:10}
L.~Scardovi, ``Clustering and synchronization in phase models with state
  dependent coupling,'' in \emph{{IEEE} Conf.\ on Decision and Control},
  Atlanta, GA, USA, Dec. 2010, pp. 627--632.

\bibitem{AF-AC-WPL:11}
A.~Franci, A.~Chaillet, and W.~Pasillas-L{\'e}pine, ``Existence and robustness
  of phase-locking in coupled {K}uramoto oscillators under mean-field
  feedback,'' \emph{Automatica}, vol.~47, no.~6, pp. 1193--1202, 2011.

\bibitem{SYH-MS:11}
S.~Y. Ha and M.~Slemrod, ``A fast-slow dynamical systems theory for the
  {K}uramoto type phase model,'' \emph{Journal of Differential Equations}, vol.
  251, no.~10, pp. 2685--2695, 2011.

\bibitem{LB-LS-ADG:09}
L.~Buzna, S.~Lozano, and A.~Diaz-Guilera, ``Synchronization in symmetric
  bipolar population networks,'' \emph{Physical Review E}, vol.~80, no.~6, p.
  66120, 2009.

\bibitem{YM-AFP:04}
Y.~Moreno and A.~F. Pacheco, ``Synchronization of {K}uramoto oscillators in
  scale-free networks,'' \emph{Europhysics Letters}, vol.~68, p. 603, 2004.

\bibitem{ACK:10}
A.~C. Kalloniatis, ``From incoherence to synchronicity in the network
  {K}uramoto model,'' \emph{Physical Review E}, vol.~82, no.~6, p. 066202,
  2010.

\bibitem{AS:09}
A.~Sarlette, ``Geometry and symmetries in coordination control,'' Ph.D.
  dissertation, University of Li\`ege, Belgium, Jan. 2009.

\bibitem{EAC-PAM-RF:10-a}
E.~A. Canale, P.~A. Monzn, and F.~Robledo, ``{The wheels: an infinite family of
  bi-connected planar synchronizing graphs},'' in \emph{IEEE Conf. Industrial
  Electronics and Applications}, Taichung, Taiwan, Jun. 2010, pp. 2204--2209.

\bibitem{EAC-PAM-RF:10-b}
E.~A. Canale, P.~Monzon, and F.~Robledo, ``On the complexity of the
  classification of synchronizing graphs,'' in \emph{Grid and Distributed
  Computing, Control and Automation}, Jeju Island, Korea, Dec. 2010, pp.
  186--195.

\bibitem{PM:06}
P.~Monz{\'o}n, ``Almost global stability of dynamical systems,'' Ph.D.
  dissertation, Universidad de la Rep{\'u}blica, Montevideo, Uruguay, Jul.
  2006.

\bibitem{JLvH-WFW:93}
J.~L. van Hemmen and W.~F. Wreszinski, ``Lyapunov function for the {K}uramoto
  model of nonlinearly coupled oscillators,'' \emph{Journal of Statistical
  Physics}, vol.~72, no.~1, pp. 145--166, 1993.

\bibitem{SJC-JJS:10b}
S.~J. Chung and J.~J. Slotine, ``{On synchronization of coupled Hopf-Kuramoto
  oscillators with phase delays},'' in \emph{{IEEE} Conf.\ on Decision and
  Control}, Atlanta, GA, USA, Dec. 2010, pp. 3181--3187.

\bibitem{GSS-UM-FA:09}
G.~S. Schmidt, U.~M{\"u}nz, and F.~Allg{\"o}wer, ``Multi-agent speed consensus
  via delayed position feedback with application to {K}uramoto oscillators,''
  in \emph{{E}uropean {C}ontrol {C}onference}, Budapest, Hungary, Aug. 2009,
  pp. 2464--2469.

\bibitem{EC-PM:09}
E.~Canale and P.~Monz{\'o}n, ``On the characterization of families of
  synchronizing graphs for {K}uramoto coupled oscillators,'' in \emph{IFAC
  Workshop on Distributed Estimation and Control in Networked Systems}, Venice,
  Italy, Sep. 2009, pp. 42--47.

\bibitem{SHS-REM:88}
S.~H. Strogatz and R.~E. Mirollo, ``Phase-locking and critical phenomena in
  lattices of coupled nonlinear oscillators with random intrinsic
  frequencies,'' \emph{Physica D: Nonlinear Phenomena}, vol.~31, no.~2, pp.
  143--168, 1988.

\bibitem{CH:10}
H.~Chiba, ``A proof of the {K}uramoto's conjecture for a bifurcation structure
  of the infinite dimensional {K}uramoto model,'' \emph{Arxiv preprint
  arXiv:1008.0249}, 2010.

\bibitem{EAM-EB-SHS-PS-TMA:09}
E.~A. Martens, E.~Barreto, S.~H. Strogatz, E.~Ott, P.~So, and T.~M. Antonsen,
  ``Exact results for the {K}uramoto model with a bimodal frequency
  distribution,'' \emph{Physical Review E}, vol.~79, no.~2, p. 26204, 2009.

\bibitem{YH-PGM-SPM-UVS:12}
H.~Yin, P.~G. Mehta, S.~P. Meyn, and U.~V. Shanbhag, ``Synchronization of
  coupled oscillators is a game,'' \emph{IEEE Transactions on Automatic
  Control}, vol.~57, no.~4, pp. 920--935, 2012.

\bibitem{DAW-SHS-MG:06}
D.~A. Wiley, S.~H. Strogatz, and M.~Girvan, ``The size of the sync basin,''
  \emph{Chaos}, vol.~16, no.~1, p. 015103, 2006.

\bibitem{YW-JFD:11b}
Y.~Wang and F.~J. Doyle, ``On influences of global and local cues on the rate
  of synchronization of oscillator networks,'' \emph{Automatica}, vol.~47,
  no.~6, pp. 1236--1242, 2011.

\bibitem{FD-MC-FB:11v-arxiv}
F.~D{\"o}rfler, M.~Chertkov, and F.~Bullo, ``Synchronization in complex
  oscillator networks and smart grids,'' Jul. 2012, available at
  \texttt{http://arxiv.org/pdf/1208.0045}.

\bibitem{EWJ-PSK:04}
E.~W. Justh and P.~S. Krishnaprasad, ``Equilibria and steering laws for planar
  formations,'' \emph{Systems \& Control Letters}, vol.~52, no.~1, pp. 25--38,
  2004.

\bibitem{SS-PV:80}
S.~Sastry and P.~Varaiya, ``Hierarchical stability and alert state steering
  control of interconnected power systems,'' \emph{IEEE Transactions on
  Circuits and Systems}, vol.~27, no.~11, pp. 1102--1112, 1980.

\bibitem{FD-FB:11d}
F.~D{\"o}rfler and F.~Bullo, ``{K}ron reduction of graphs with applications to
  electrical networks,'' \emph{IEEE Transactions on Circuits and Systems}, Nov.
  2011, to appear.

\bibitem{VF-SR-EC-EM-VR:09}
V.~Fioriti, S.~Ruzzante, E.~Castorini, E.~Marchei, and V.~Rosato, ``Stability
  of a distributed generation network using the {K}uramoto models,'' in
  \emph{Critical Information Infrastructure Security}, ser. Lecture Notes in
  Computer Science.\hskip 1em plus 0.5em minus 0.4em\relax Springer, 2009, pp.
  14--23.

\bibitem{GF-AHN-NFP:08}
G.~Filatrella, A.~H. Nielsen, and N.~F. Pedersen, ``Analysis of a power grid
  using a {K}uramoto-like model,'' \emph{The European Physical Journal~B},
  vol.~61, no.~4, pp. 485--491, 2008.

\bibitem{MR-AS-MT-DW:12}
M.~Rohden, A.~Sorge, M.~Timme, and D.~Witthaut, ``Self-organized
  synchronization in decentralized power grids,'' \emph{Physical Review
  Letters}, vol. 109, no.~6, p. 064101, 2012.

\bibitem{DS-UR-BS-MP:01}
D.~Subbarao, R.~Uma, B.~Saha, and M.~V.~R. Phanendra, ``Self-organization on a
  power system,'' \emph{IEEE Power Engineering Review}, vol.~21, no.~12, pp.
  59--61, 2001.

\bibitem{DJH-GC:06}
D.~J. Hill and G.~Chen, ``Power systems as dynamic networks,'' in \emph{IEEE
  Int. Symposium on Circuits and Systems}, Kos, Greece, May 2006, pp. 722--725.

\bibitem{FCH-EMI:97}
F.~C. Hoppensteadt and E.~M. Izhikevich, \emph{Weakly connected neural
  networks}.\hskip 1em plus 0.5em minus 0.4em\relax Springer, 1997, vol. 126.

\bibitem{EMI:07}
E.~M. Izhikevich, \emph{Dynamical Systems in Neuroscience: The Geometry of
  Excitability and Bursting}.\hskip 1em plus 0.5em minus 0.4em\relax MIT Press,
  2007.

\bibitem{EMI-YK:06}
E.~M. Izhikevich and Y.~Kuramoto, ``Weakly coupled oscillators,''
  \emph{Encyclopedia of Mathematical Physics}, vol.~5, p. 448, 2006.

\bibitem{GBE-NK:84}
G.~B. Ermentrout and N.~Kopell, ``Frequency plateaus in a chain of weakly
  coupled oscillators, {I.}'' \emph{SIAM journal on Mathematical Analysis},
  vol.~15, no.~2, pp. 215--237, 1984.

\bibitem{FFW-SK:80}
F.~F. Wu and S.~Kumagai, \emph{Limits on Power Injections for Power Flow
  Equations to Have Secure Solutions}.\hskip 1em plus 0.5em minus 0.4em\relax
  Electronics Research Laboratory, College of Engineering, University of
  California, 1980.

\bibitem{TN-AEM-YCL-FCH:03}
T.~Nishikawa, A.~E. Motter, Y.~C. Lai, and F.~C. Hoppensteadt, ``Heterogeneity
  in oscillator networks: {A}re smaller worlds easier to synchronize?''
  \emph{Physical Review Letters}, vol.~91, no.~1, p. 14101, 2003.

\bibitem{FW-SK:82}
F.~Wu and S.~Kumagai, ``Steady-state security regions of power systems,''
  \emph{IEEE Transactions on Circuits and Systems}, vol.~29, no.~11, pp.
  703--711, 1982.

\bibitem{GK-MBH-KEB-MJB-BK-DA:06}
G.~Korniss, M.~B. Hastings, K.~E. Bassler, M.~J. Berryman, B.~Kozma, and
  D.~Abbott, ``Scaling in small-world resistor networks,'' \emph{Physics
  Letters A}, vol. 350, no. 5-6, pp. 324--330, 2006.

\bibitem{LMP-TLC:98}
L.~M. Pecora and T.~L. Carroll, ``Master stability functions for synchronized
  coupled systems,'' \emph{Physical Review Letters}, vol.~80, no.~10, pp.
  2109--2112, 1998.

\bibitem{JGG-YM-AA:07}
J.~G{\'o}mez-Gardenes, Y.~Moreno, and A.~Arenas, ``Paths to synchronization on
  complex networks,'' \emph{Physical Review Letters}, vol.~98, no.~3, p. 34101,
  2007.

\bibitem{CJT-OJMS:72a}
C.~J. Tavora and O.~J.~M. Smith, ``Stability analysis of power systems,''
  \emph{IEEE Transactions on Power Apparatus and Systems}, vol.~91, no.~3, pp.
  1138--1144, 1972.

\bibitem{AA-SS-VP:81}
A.~Araposthatis, S.~Sastry, and P.~Varaiya, ``Analysis of power-flow
  equation,'' \emph{International Journal of Electrical Power \& Energy
  Systems}, vol.~3, no.~3, pp. 115--126, 1981.

\bibitem{CJT-OJMS:72b}
C.~J. Tavora and O.~J.~M. Smith, ``Equilibrium analysis of power systems,''
  \emph{IEEE Transactions on Power Apparatus and Systems}, vol.~91, no.~3, pp.
  1131--1137, 1972.

\bibitem{MI:92}
M.~Ili\'c, ``Network theoretic conditions for existence and uniqueness of
  steady state solutions to electric power circuits,'' in \emph{IEEE
  International Symposium on Circuits and Systems}, San Diego, CA, USA, May
  1992, pp. 2821--2828.

\bibitem{KSC-DJH:86}
K.~S. Chandrashekhar and D.~J. Hill, ``Cutset stability criterion for power
  systems using a structure-preserving model,'' \emph{International Journal of
  Electrical Power \& Energy Systems}, vol.~8, no.~3, pp. 146--157, 1986.

\bibitem{RS-AS-PR:10}
R.~Sepulchre, A.~Sarlette, and P.~Rouchon, ``Consensus in non-commutative
  spaces,'' in \emph{{IEEE} Conf.\ on Decision and Control}, Atlanta, GA, USA,
  Dec. 2010, pp. 6596--6601.

\bibitem{HKK:02}
H.~K. Khalil, \emph{Nonlinear Systems}, 3rd~ed.\hskip 1em plus 0.5em minus
  0.4em\relax Prentice Hall, 2002.

\bibitem{LM:08-arxiv}
L.~Moreau, ``Stability of continuous-time distributed consensus algorithms,''
  Feb. 2008, available at \texttt{http://arxiv.org/abs/math/0409010v1}.

\bibitem{DA:02}
D.~Angeli, ``A {L}yapunov approach to incremental stability properties,''
  \emph{IEEE Transactions on Automatic Control}, vol.~47, no.~3, pp. 410--421,
  2002.

\bibitem{AM:11}
A.~Mauroy, ``On the dichotomic collective behaviors of large populations of
  pulse-coupled firing oscillators,'' Ph.D. dissertation, University of
  Li\`{e}ge, Belgium, 2011.

\bibitem{EHS:94}
E.~H. Spanier, \emph{Algebraic Topology}.\hskip 1em plus 0.5em minus
  0.4em\relax Springer, 1994.

\end{thebibliography}


%

\end{document}